\documentclass[12pt]{amsart}
\usepackage{amscd,amssymb}
\usepackage[arrow,matrix]{xy}
\usepackage[colorlinks,plainpages,backref,urlcolor=blue]{hyperref}

\newtheorem{theorem}[subsection]{Theorem}
\newtheorem{lemma}[subsection]{Lemma}
\newtheorem{prop}[subsection]{Proposition}
\newtheorem{corollary}[subsection]{Corollary}

\newtheorem{thm}{Theorem}

\theoremstyle{definition}
\newtheorem{remark}[subsection]{Remark}
\newtheorem{definition}[subsection]{Definition}
\newtheorem{example}[subsection]{Example}
\newtheorem*{ack}{Acknowledgments}

\numberwithin{equation}{section}
\newcommand{\F}{{\mathcal F}}
\newcommand{\A}{{\mathcal A}}
\newcommand{\OO}{{\mathcal O}}

\newcommand{\W}{{\mathcal W}}
\newcommand{\CC}{{\mathcal C}}

\newcommand{\V}{{\mathcal V}}
\newcommand{\R}{{\mathcal R}}
\newcommand{\TT}{{\mathcal T}}

\newcommand{\h}{\mathfrak {H}}
\newcommand{\m}{\mathfrak {m}}
\newcommand{\wh}{\widehat{\mathfrak {H}}}

\newcommand{\oM}{\overline{M}}
\newcommand{\oC}{\overline{C}}

\newcommand{\sV}{{\sf V}}
\newcommand{\sE}{{\sf E}}
\newcommand{\sW}{{\sf W}}
\newcommand{\kk}{\kappa}

\newcommand{\T}{\mathbb{T}}
\newcommand{\Z}{\mathbb{Z}}
\newcommand{\N}{\mathbb{N}}
\newcommand{\Q}{\mathbb{Q}}
\newcommand{\RR}{\mathbb{R}}
\newcommand{\C}{\mathbb{C}}
\newcommand{\K}{\Bbbk}
\newcommand{\PP}{\mathbb{P}}
\newcommand{\FF}{\mathbb{F}}
\newcommand{\bF}{\mathbb{F}}
\newcommand{\bL}{\mathbb{L}}
\newcommand{\eexp}{\mathsf{exp}}

\DeclareMathOperator{\Hom}{Hom}
\DeclareMathOperator{\rank}{rank}
\DeclareMathOperator{\im}{im}
\DeclareMathOperator{\init}{in}
\DeclareMathOperator{\coker}{coker}
\DeclareMathOperator{\spn}{span}
\DeclareMathOperator{\id}{id}

\DeclareMathOperator{\tor}{Tor}
\DeclareMathOperator{\codim}{codim}

\DeclareMathOperator{\gr}{gr}
\DeclareMathOperator{\ab}{ab}
\DeclareMathOperator{\ad}{ad}
\DeclareMathOperator{\odd}{odd}
\DeclareMathOperator{\Zero}{Zero}
\DeclareMathOperator{\Sym}{Sym}
\DeclareMathOperator{\Sing}{Sing}

\newcommand{\surj}{\twoheadrightarrow}

\newcommand{\llangle}{{\langle\!\langle}}
\newcommand{\rrangle}{{\rangle\!\rangle}}
\newcommand{\abs}[1]{\left| #1 \right|}

\newenvironment{romenum}
{

\begin{enumerate}}{\end{enumerate}}

\begin{document}

\title[Formality, Alexander invariants, and a question of Serre]%
{Formality, Alexander invariants, and \\[3pt] a question of Serre}

\author[A. Dimca]{Alexandru Dimca}
\address{  Laboratoire J.A.~Dieudonn\'{e}, UMR du CNRS 6621, 
                 Universit\'{e} de Nice-Sophia-Antipolis,
                 Parc Valrose,
                 06108 Nice Cedex 02,
                 France}
\email
{dimca@math.unice.fr}

\author[\c{S}. Papadima ]{\c{S}tefan Papadima$^{1}$}
\address{Inst. of Math. ``Simion Stoilow",
P.O. Box 1-764,
RO-014700 Bucharest, Romania}
\email
{Stefan.Papadima@imar.ro}

\author[A.~I.~Suciu]{Alexander~I.~Suciu$^{2}$}
\address{Department of Mathematics,
Northeastern University,
Boston, MA 02115, USA}
\email{a.suciu@neu.edu}
\urladdr{http://www.math.neu.edu/\~{}suciu}

\thanks{$^1$Partially supported by the CEEX Programme 
of the Romanian Ministry of Education and Research,
contract 2-CEx 06-11-20/2006}
\thanks{$^2$Partially supported by NSF grant DMS-0311142}



\begin{abstract}
We elucidate the key role played by formality in the 
theory of characteristic and resonance varieties. 
We show that the $I$-adic completion of the Alexander 
invariant of a $1$-formal group $G$ is determined solely 
by the cup-product map in low degrees. It follows that the 
germs at the origin of the characteristic and resonance 
varieties of $G$ are analytically isomorphic; in particular, 
the tangent cone  to $V_k(G)$ at $1$ equals $R_k(G)$.  
This provides new obstructions to $1$-formality.   
A detailed analysis of the irreducible components of
the tangent cone at $1$ to the first characteristic 
variety yields powerful obstructions to realizing a 
finitely presented group as the fundamental group 
of a smooth, complex quasi-projective algebraic variety.  
This sheds new light on a classical problem of J.-P.~Serre. 
Applications to arrangements, configuration spaces, 
coproducts of groups, and Artin groups are given.
\end{abstract}

\maketitle

\tableofcontents

\section{Introduction} 
\label{sec=intro}

A recurring theme in algebraic topology and  
geometry is the extent to which the cohomology of a 
space reflects the underlying topological or geometric 
properties of that space.  In this paper, we focus on the 
degree-one cohomology jumping loci of the fundamental 
group $G$ of a connected, finite-type CW-complex $M$: 
the characteristic varieties $\V_k(G)$, and the resonance 
varieties $\R_k(G)$.  Our goal is to relate these two 
sets of varieties, and to better understand their structural 
properties, under certain naturally defined conditions 
on $M$.  In turn, this analysis yields new and powerful 
obstructions for a finitely presented group $G$ to be 
$1$-formal, or realizable as the fundamental 
group of a smooth, complex quasi-projective variety.

\subsection{Cohomology jumping loci}
\label{subsec=cohojumps}
Let $\T_G=\Hom(G,\C^*)$ be the character variety of 
$G=\pi_1(M)$.  The {\em characteristic varieties}\/ of $M$ 
are the jumping loci for the cohomology of $M$, with 
coefficients in rank~$1$ local systems:
\begin{equation} 
\label{eq:charvarx}
\V^i_k(M)=\{\rho \in \T_G \mid \dim H^i(M, {}_\rho\C)\ge k\}.
\end{equation}

These varieties emerged from the work of S.~Novikov \cite{N}
on Morse theory for closed $1$-forms on manifolds.  
Foundational results on the structure of the cohomology 
support loci for local systems on quasi-projective algebraic 
varieties were obtained by  Beauville \cite{Beau1, Beau2}, 
Green and Lazarsfeld \cite{GL}, Simpson \cite{Si92, Sim}, 
Arapura \cite{A}, and Campana \cite{Cam}.

Let $H^{\bullet} (M,\C)$ be the cohomology algebra 
of $M$.  Right-multiplication by an element $z\in H^1(M,\C)$ 
yields a cochain complex $(H^{\bullet}(M,\C),\mu_z)$.  
The {\em resonance varieties}\/ of $M$ are the jumping 
loci for the homology of this complex:
\begin{equation} 
\label{eq:resvarx}
\R^i_k(M)=\{z \in H^1(M,\C) \mid 
\dim H^i(H^{\bullet}(M,\C),\mu_z) \ge  k\}.
\end{equation}

These homogeneous algebraic subvarieties of 
$H^1(M,\C)=\Hom(G,\C)$ 
were first defined by Falk \cite{F} in the case when $M$ is the 
complement of a complex hyperplane arrangement; in this 
setting, a purely combinatorial description of $\R^1_k(M)$ 
was given by Falk \cite{F}, Libgober--Yuzvinsky \cite{LY},
and Falk--Yuzvinsky \cite{FY}.

We consider here only the cohomology jumping loci in 
degree $i=1$.  These loci depend exclusively on  
$G=\pi_1(M)$, so we write $\V_k(G)=\V^1_k(M)$ and 
$\R_k(G)=\R^1_k(M)$.  The higher degree jumping loci  
will be treated in a forthcoming paper.

\subsection{The tangent cone theorem}
\label{subsec=tangent}

The key topological property that allows us to relate the 
characteristic and resonance varieties of a space $M$ 
is {\em formality}, in the sense of D.~Sullivan \cite{S}.  
Since we deal solely with the fundamental group 
$G=\pi_1(M)$, we only need $G$ to be {\em $1$-formal}. 
This property requires that $E_G$, the Malcev 
Lie algebra of $G$ (over $\C$), be isomorphic, 
as a filtered Lie algebra, to the holonomy Lie algebra 
$\h(G)=\bL/\langle \im \partial_G\rangle$, completed 
with respect to bracket length filtration.  
Here $\bL$ denotes the free Lie algebra on $H_1(G,\C)$, 
while $\partial_G$ denotes the dual of the cup-product map 
$\cup _G\colon \bigwedge ^2H^1(G,\C) \to H^2(G,\C)$. 

A group $G$ is $1$-formal if and only if $E_G$ can  
be presented with quadratic relations only; see 
Section \ref{sec=holo} for details. Many interesting 
groups fall in this class: knot groups (\cite{S}), 
certain pure braid groups of Riemann surfaces 
(Bezrukavnikov \cite{B}, Hain \cite{Hai}),  fundamental 
groups of compact K\"{a}hler manifolds (Deligne, Griffiths, 
Morgan, and Sullivan \cite{DGMS}), fundamental groups 
of complements of hypersurfaces in $\C\PP^n$ 
(Kohno \cite{K}), and finite-type Artin groups 
(Kapovich--Millson \cite{KM}) are all $1$-formal. 

\begin{thm} 
\label{thm=tcfintro} 
Let $G$ be a $1$-formal group. For each 
$0\le k \le b_1(G)$, the exponential map 
$\exp\colon \Hom(G,\C) \to \Hom(G,\C^*)$ 
restricts to an isomorphism of analytic germs,  
\[
\exp\colon (\R_k(G),0)  \xrightarrow{\,\simeq\,} (\V_k(G),1).
\]
In particular, the tangent cone at $1$ to $\V_k(G)$
equals $\R_k(G)$:
\[
TC_1(\V_k(G))=\R_k(G) \, .
\]
\end{thm}

Essential ingredients in the proof are two 
modules associated to a finitely presented 
group $G$:  the {\em Alexander invariant}, 
$B_G=G'/G''$, viewed as a module over the group 
ring $\Z G_{\ab}$, and its ``infinitesimal" version, 
$B_{\h(G)}=\coker (\delta_3 + \partial_G)$, 
viewed as a module over the polynomial ring 
$\C[X]=\Sym(G_{\ab}\otimes \C)$, where $\delta_j$ 
denotes the $j$-th Koszul differential. The 
varieties defined by the Fitting ideals of these 
two modules coincide, away from the origin, 
with $\V_k(G)$ and $\R_k(G)$, respectively.

Under our $1$-formality assumption on $G$, we deduce 
Theorem \ref{thm=tcfintro} from the fact that the exponential 
map induces a filtered isomorphism between the $I$-adic 
completion of $B_G \otimes \C$ and the $(X)$-adic 
completion of $B_{\h(G)}$. We establish this 
isomorphism using Malcev Lie algebra tools. Related  
techniques have been used previously in low-dimensional 
topology \cite{Hai, P2}, algebraic geometry \cite{M, ABC},  
and group theory \cite{MP, PS2}. 

Theorem \ref{thm=tcfintro} sharpens and extends several  
results from \cite{ESV, STV, CS2}, which only 
apply to the case when $G$ is the fundamental group 
of the complement of a complex hyperplane arrangement.  
The point is that only a {\em topological}\/ property 
($1$-formality) is needed for the conclusion to hold.  
Further information in the case of hypersurface arrangements 
can be found in \cite{Li01, DM}. 

For an arbitrary finitely presented group $G$, the tangent cone 
to $\V_k(G)$ at the origin is contained in $\R_k(G)$, see 
Libgober \cite{Li}.  Yet the inclusion may well be strict.  
In fact, as noted in Example  \ref{ex:notinj} and 
Remark \ref{rem:lib}, the tangent cone formula 
from Theorem \ref{thm=tcfintro} fails even for 
fundamental groups of smooth, quasi-projective 
varieties.  

Theorem \ref{thm=tcfintro} provides a new, and quite 
powerful obstruction to $1$-formality of groups---and thus, 
to formality of spaces.  As illustrated in Example \ref{ex=nonformal}, 
this obstruction complements, and in some cases strengthens, 
classical obstructions to ($1$-) formality, due to Sullivan, 
such as the existence of an isomorphism 
$\gr(G) \otimes \C \cong \h(G)$.   

\subsection{Serre's question}
\label{subsec=serre}

As is well-known, every finitely presented group $G$ is the 
fundamental group of a smooth, compact, connected 
$4$-dimensional manifold $M$.  Requiring a complex structure 
on $M$ is no more restrictive, as long as one is willing to go up 
in dimension; see Taubes \cite{Ta}. Requiring that $M$ be a 
compact K\"{a}hler manifold, though, puts extremely strong 
restrictions on what $G=\pi_1(M)$ can be.  
We refer to \cite{ABC} for a comprehensive survey 
of K\"{a}hler groups.   

J.-P. Serre asked the following question: 
which finitely presented groups can be realized as 
fundamental groups of smooth, connected, 
quasi-projective, complex algebraic varieties? 
Following Catanese \cite{Cat03}, 
we shall call a group $G$ arising in this fashion a 
{\em quasi-projective}\/ group.

In this context, one may also consider the  
larger class of quasi-compact K\"{a}hler manifolds, 
of the form $M= \oM \setminus D$, where 
$\oM$ is compact K\"{a}hler and $D$ is a normal 
crossing divisor.  If $G=\pi_1(M)$ with $M$ as above, 
we say $G$ is a {\em quasi-K\"{a}hler} group.

The first obstruction to quasi-projectivity was given 
by J. Morgan:   If $G$ is a quasi-projective group, then 
$E_G=\widehat{\bL}/J$, where $\bL$ is a free 
Lie algebra with generators in degrees $1$ and $2$,  
and $J$ is a closed Lie ideal, generated in degrees 
$2$, $3$ and $4$; see \cite[Corollary 10.3]{M}.  
By refining Morgan's techniques, Kapovich and Millson 
obtained analogous quasi-homogeneity restrictions, 
on certain singularities of representation varieties of $G$ 
into reductive algebraic Lie groups; see \cite[Theorem 1.13]{KM}.  
Another obstruction is due to Arapura:  If $G$ is  
quasi-K\"{a}hler, then the characteristic variety $\V_1(G)$ 
is a union of (possibly translated) subtori of $\T_G$; 
see \cite[p.~564]{A}. 

If the group $G$ is $1$-formal, then $E_G=\widehat{\bL}/J$, 
with $\bL$ generated in degree $1$ and $J$ generated in 
degree $2$; thus, $G$ verifies Morgan's test.  It is therefore 
natural to explore the relationship between $1$-formality 
and quasi-projectivity.  (In contrast with the K\"{a}hler case, 
it is known from \cite{M, KM} that these two properties are 
independent.) Another motivation for our investigation comes 
from the study of fundamental groups of complements of plane 
algebraic curves. This class of $1$-formal, quasi-projective 
groups contains, among others, the celebrated Stallings group;  
see \cite{PS06b}.

\subsection{Position and resonance obstructions}
\label{subsec=resobs}  

Our main contribution to Serre's problem is Theorem 
\ref{thm=posobstr} below, which provides a new type of 
restriction on fundamental groups of smooth, quasi-projective 
complex algebraic varieties. In the presence of $1$-formality, 
this restriction is expressed entirely in terms of a 
classical invariant, namely the cup-product map 
$\cup _G\colon \bigwedge ^2H^1(G,\C) \to H^2(G,\C)$. 

\begin{thm} 
\label{thm=posobstr} 
Let $M$ be a  connected, quasi-compact K\"{a}hler manifold. 
Set $G=\pi_1(M)$ and let $\{ \V^{\alpha} \}$ be the 
irreducible components of $\V_1(G)$ containing $1$. Denote 
by $\TT^{\alpha}$ the tangent space at $1$ of $\V^{\alpha}$.
Then the following hold.

\begin{enumerate}
\item  \label{a1}  
Every positive-dimensional tangent space $\TT^{\alpha}$ 
is a $p$-isotropic linear subspace of $H^1(G, \C)$, 
of dimension at least $2p+2$, for some 
$p=p(\alpha) \in \{0,1\}$. 

\item  \label{a2} 
If $\alpha \ne \beta $, then $\TT^{\alpha} \bigcap \TT^{\beta}=\{0\}$.
\end{enumerate}
Assume further that $G$ is $1$-formal.  Let $\{ \R^{\alpha} \}$ 
be the irreducible components of $\R_1(G)$.  Then 
the following hold.
\begin{enumerate}
\setcounter{enumi}{2}
\item \label{a0}
The collection $\{ \TT^{\alpha}\}$ coincides with the 
collection $\{ \R^{\alpha}\}$.

\item  \label{a3}
For $1\leq k\leq b_1(G)$, we have 
$\R_k(G)=\{0\} \cup \bigcup\nolimits_\alpha \R^{\alpha}$, 
where the union is over all components  $\R^{\alpha}$ such that 
$\dim \R^{\alpha} >k+p(\alpha)$.  

\item \label{a4}
The group $G$ has a free quotient of rank at least $2$ 
if and only if $\R_1(G)$ strictly contains $\{ 0\}$.
\end{enumerate}

\end{thm}

Here, we say that a non-zero subspace $V\subseteq  H^1(G, \C)$ 
is $0$- (respectively, $1$-) {\em isotropic}\/ if the restriction 
of the cup-product map, 
$\cup_G\colon \bigwedge^2 V\to \cup_G (\bigwedge^2 V)$, 
is equivalent to 
$\cup_C \colon \bigwedge^2 H^1(C, \C)\to H^2(C, \C)$,
where $C$ is a non-compact (respectively, compact) smooth, 
connected complex curve. See \ref{def=position} for a more 
precise definition.

In this paper, we consider only components of $\V_1(G)$ 
containing $1$. For a detailed analysis of translated components, 
we refer to \cite{D2} and \cite{DPS-codone}. 

The proof of Theorem \ref{thm=posobstr} is given in 
Section \ref{sec=posobs}, and relies on results of Arapura \cite{A} 
on quasi-K\"{a}hler groups. Part \eqref{a1} is an easy consequence 
of \cite[Proposition V.1.7]{A}. Part \eqref{a2} is a new viewpoint, 
developed in \cite{DPS-codone} to obtain a completely new type 
of obstruction.  Part \eqref{a0} follows from our tangent cone formula.

For an arrangement complement $M$, Parts \eqref{a2} and 
\eqref{a3} of the above theorem were proved by Libgober 
and Yuzvinsky in \cite{LY}, using completely different methods. 

The ``position" obstructions \eqref{a1} and \eqref{a2} in
Theorem \ref{thm=posobstr} can be viewed as a 
much strengthened form of Arapura's theorem:
they give information on 
how the components of $\V_1(G)$ passing through 
the origin intersect at that point, and how their tangent 
spaces at $1$ are situated with respect to the 
cup-product map $\cup_G$.   

We may also
view conditions \eqref{a1}--\eqref{a4} as a set of ``resonance" 
obstructions for a $1$-formal group to be quasi-K\"{a}hler, 
or for a quasi-K\"{a}hler group to be $1$-formal.  
Since the class of homotopy types of compact K\"{a}hler 
manifolds is strictly larger than the class of homotopy types 
of smooth projective varieties (see Voisin \cite{V}), one 
may wonder whether the class of quasi-K\"{a}hler 
groups is also strictly larger than the class of 
quasi-projective groups.

\subsection{Applications}
\label{subsec=apps}
In the last four sections, we illustrate the efficiency of our  
obstructions to $1$-formality and quasi-K\"{a}hlerianity with 
several classes of examples.  

We start in Section \ref{sec=wp} with wedges and products 
of spaces. Our analysis of resonance varieties of wedges,  
together with Theorem \ref{thm=posobstr},  shows that 
$1$-formality and quasi-K\"{a}hlerianity behave quite differently 
with respect to the coproduct operation for groups. 
Indeed, if $G_1$ and $G_2$ are $1$-formal, 
then $G_1*G_2$ is also $1$-formal; but, if in addition 
both factors are non-trivial, presented by 
commutator relators only, and one of them is non-free,
then $G_1*G_2$ is not quasi-K\"{a}hler. 
As a consequence of the position obstruction from Theorem
\ref{thm=posobstr}\eqref{a1}, we also show that the
quasi-K\"{a}hlerianity of $G_1*G_2$, where the groups $G_i$
are assumed only finitely presented with infinite abelianization,
implies the vanishing of $\cup_{G_1}$ and $\cup_{G_2}$.

When it comes to resonance varieties, real subspace 
arrangements offer a stark contrast to complex 
hyperplane arrangements, cf.~\cite{MS1, MS2}.  
If $M$ is the complement of an arrangement of transverse 
planes through the origin of $\RR^4$, then $G=\pi_1(M)$ 
passes Sullivan's $\gr$-test.  Yet, as we note in 
Section \ref{sec=realarr}, the group $G$ may fail 
the tangent cone formula from Theorem \ref{thm=tcfintro}, 
and thus be non-$1$-formal; or, $G$ may be $1$-formal, 
but fail tests \eqref{a1}, \eqref{a2}, \eqref{a3} from 
Theorem \ref{thm=posobstr}, and thus be 
non-quasi-K\"{a}hler.

In Section \ref{sec=conf}, we apply our techniques to the 
configuration spaces $M_{g,n}$ of $n$ ordered points 
on a closed Riemann surface of genus $g$.  Clearly, 
the surface pure braid groups $P_{g,n}=\pi_1(M_{g,n})$ 
are quasi-projective.  On the other hand, if $n\ge 3$, 
the variety $\R_1(P_{1,n})$ is irreducible and non-linear. 
Theorem \ref{thm=posobstr}\eqref{a0} shows that 
$P_{1,n}$ is not $1$-formal, thereby verifying a 
statement of Bezrukavnikov \cite{B}.  

We conclude in Section \ref{sec=artinalg} with a study 
of Artin groups associated to finite, labeled graphs, 
from the perspective of their cohomology jumping loci.  
As shown in \cite{KM}, all Artin groups are $1$-formal; 
thus, they satisfy Morgan's homogeneity test. Moreover, 
as we show in Proposition \ref{prop=vart} (building on 
work from \cite{PS1}), the first characteristic varieties  
of right-angled Artin groups are unions of coordinate 
subtori; thus, all such groups pass Arapura's $\V_1$-test.

In \cite[Theorem 1.1]{KM}, Kapovich and Millson establish, 
by a fairly involved argument, the existence of infinitely many 
Artin groups that are not realizable by smooth, quasi-projective 
varieties. Using the isotropicity obstruction from Theorem 
\ref{thm=posobstr}, we show that a right-angled Artin group 
$G_{\Gamma}$ is quasi-K\"{a}hler  if and only if $\Gamma$ 
is a complete, multi-partite graph, in which case $G_{\Gamma}$ 
is actually quasi-projective.  This result provides a complete---and 
quite satisfying---answer to Serre's problem within this class of 
groups.  In the process, we take a first step towards 
solving the problem for all Artin groups, by answering it at 
the level of associated Malcev Lie algebras. We also 
determine precisely which right-angled Artin groups 
are K\"{a}hler.

The approach to Serre's problem taken in this paper---based 
on the obstructions from Theorem \ref{thm=posobstr}---has 
led to complete answers for several other classes of groups:
\begin{itemize}
\item In \cite{DPS-bb}, we classify the quasi-K\"{a}hler 
groups within the class of groups introduced by Bestvina 
and Brady in \cite{BB}.\\[-5pt]
\item In \cite{DPS-codone}, we determine precisely which 
fundamental groups of boundary manifolds of line 
arrangements in $\C\PP^2$ are quasi-projective groups.\\[-5pt]
\item In \cite{DS}, we decide which $3$-manifold groups 
are K\"{a}hler groups, thus answering a question raised by 
S.~Donaldson and W.~Goldman in 1989.\\[-5pt]
\item In \cite{P07}, we show that the fundamental groups of 
a certain natural class of graph links in $\Z$-homology spheres 
are quasi-projective if and only if the corresponding links are 
Seifert links.
\end{itemize}

The obstructions from Theorem \ref{thm=posobstr} are 
complemented by new methods of constructing 
interesting (quasi-)projective groups.  
These methods, developed in \cite{DPS-bb} and 
\cite{DPS-complexbb}, lead to a negative answer to the 
following question posed by J.~Koll\'{a}r in \cite{Ko}:  
Is every projective group commensurable (up to finite 
kernels) with a group admitting a quasi-projective 
variety as classifying space?  

\section{Holonomy Lie algebra, Malcev completion, 
and $1$-formality} 
\label{sec=holo}

Given a (finitely presented) group $G$, there are several 
Lie algebras attached to it:  the associated graded Lie 
algebra $\gr^*(G)$, the holonomy Lie algebra $\h(G)$, 
and the Malcev Lie algebra $E_G$.   
In this section, we review the constructions of these 
Lie algebras, and the related notion of $1$-formality, 
which will be crucial in what follows.

\subsection{ Holonomy Lie algebras} 
\label{subsecz=holo}
We start by recalling the definition of the holonomy 
Lie algebra, due to K.-T.~Chen \cite{C}. 

Let $M$ be a connected CW-complex with finite 
$2$-skeleton.  Let $\K$ be a field of characteristic $0$.  
Denote by $\bL^*(H_1M)$ the free Lie algebra on 
$H_1M=H_1(M,\K)$, graded by bracket length, and 
use the Lie  bracket to identify $H_1M \wedge H_1M$ 
with  $\bL^2(H_1M)$. Set
\begin{equation} 
\label{eq=holm}
\h (M):=\bL^*(H_1M)/\langle \im (\partial_M\colon 
H_2M  \to  \bL^2(H_1M))\rangle , 
\end{equation}
where $\partial_M$ is the dual of the cup-product map 
$\cup_M\colon H^1(M,\K) \wedge H^1(M,\K) \to H^2(M,\K)$ and 
$\langle (\cdot) \rangle$ denotes the Lie ideal spanned by 
$(\cdot)$. Note that $\h (M)$ is a quadratic Lie algebra, in that 
it has a presentation with generators in degree $1$ and 
relations in degree $2$ only. We call $\h (M)$ the 
{\em holonomy Lie algebra}\/ of $M$ (over the field $\K$).

Now let $G$ be a group admitting a finite presentation. 
Choose an Eilenberg-MacLane space $K(G,1)$ with 
finite $2$-skeleton, and define
\begin{equation} 
\label{eq=holg}
\h (G):=\h (K(G,1)) .
\end{equation}
If $M$ is a CW-complex as above, with $G=\pi_1(M)$, and 
if $f\colon M \to K(G,1)$ is a classifying map, then $f$ induces 
an isomorphism on $H_1$ and an epimorphism on $H_2$. 
This implies that
\begin{equation} 
\label{eq=holeq}
\h (G)=\h (M) .
\end{equation}

\subsection{Malcev Lie algebras} 
\label{subsec=malcev}
Next, we recall some useful facts about the functorial Malcev 
completion of a group, following Quillen \cite[Appendix A]{Q}.  

A {\em Malcev Lie algebra}\/ is a $\K$-Lie algebra $E$, endowed 
with a decreasing, complete filtration
\[
E=F_1E \supset \cdots \supset F_qE \supset  F_{q+1}E \supset \cdots ,
\]
by $\K$-vector subspaces  satisfying $[F_rE,F_sE] \subset  F_{r+s}E$ 
for all $r,s \ge 1$, and with the property that the associated graded 
Lie algebra, $\gr_F^*(E)=\bigoplus_{q\ge 1} F_qE /F_{q+1}E$, is 
generated by $\gr_F^1(E)$.  By completeness of the filtration, the 
Campbell-Hausdorff formula
\begin{equation} 
\label{eq=ch}
e \cdot f=e+f+\tfrac{1}{2}[e,f]+ 
\tfrac{1}{12}([e,[e,f]]+[f,[f,e]]) + \cdots 
\end{equation}
endows the underlying set of $E$ with 
a group structure, to be denoted by $\eexp(E)$.

For a group $G$, denote by $\wh (G)$ the completion of the 
holonomy Lie algebra with respect to the degree filtration.  
It is readily checked that $\wh (G)$, together with the canonical 
filtration of the completion, is a Malcev Lie algebra with 
$\gr_F^*(\wh (G))=\h^* (G)$.

In \cite{Q}, Quillen associates to $G$, in a functorial 
way, a Malcev Lie algebra $E_G$ and a group homomorphism 
$\kk _G\colon G \to \eexp(E_G)$.  The key property of the Malcev 
completion $(E_G,\kappa_G)$ is that  $\kk _G$ induces an 
isomorphism of graded $\K$-Lie algebras
\begin{equation} 
\label{eq=griso}
\gr^*( \kk _G)\colon \gr^*(G) \otimes \K  
\xrightarrow{\,\simeq\,} \gr^*_F(E_G) .
\end{equation}
Here, $\gr^*(G)=\bigoplus_{q\ge 1} \Gamma_qG /\Gamma_{q+1} G$ 
is the graded Lie algebra associated to the lower central series,
$G =\Gamma _1G \supset \cdots \supset \Gamma _qG  
\supset \Gamma _{q+1}G \supset \cdots$, defined inductively 
by setting $ \Gamma _{q+1}G=(G, \Gamma _qG)$, where 
$(A,B)$ denotes the subgroup generated by all commutators 
$(g,h)=ghg^{-1}h^{-1}$ with $g\in A$ and $h\in B$, and with the 
Lie bracket on $\gr^*(G)$ induced by the commutator map 
$(\:,\:)\colon G\times G\to G$.  

\subsection{Formal spaces and $1$-formal groups} 
\label{subsec=formal}

The important notion of formality of a space was introduced by 
D.~Sullivan in \cite{S}.  Let $M$ be a connected CW-complex 
with finitely many cells in each dimension. Roughly speaking, 
$M$ is formal if the rational cohomology algebra of $M$ 
determines the homotopy type of $M$ modulo torsion.  
Many interesting spaces are formal: compact K\"{a}hler 
manifolds \cite{DGMS}, homogeneous spaces of compact 
connected Lie groups with equal ranks \cite{S}; 
products and wedges of formal spaces are again formal.

\begin{definition} 
\label{def=1formal} 
A finitely presented group $G$ is {\em $1$-formal}\/  
if its Malcev Lie algebra, $E_G$, is isomorphic to the 
completion of its holonomy Lie algebra, $\wh (G)$, 
as filtered Lie algebras.
\end{definition}

A fundamental result of Sullivan \cite{S} states that 
$\pi _1(M)$ is $1$-formal whenever $M$ is formal. The converse 
is not true in general, though it holds under certain additional 
assumptions, see \cite{PS1}.  Here are some motivating 
examples. 

\begin{example} 
\label{ex:w2formal}
If $M$ is obtained from a smooth, complex projective 
variety $\oM$ by deleting a subvariety 
$D \subset \oM$ with $\codim D \geq 2$, 
then $\pi_1(M)=\pi_1(\oM)$. Hence, $\pi_1(M)$ 
is $1$-formal, by \cite{DGMS}.
\end{example}

\begin{example} 
\label{ex:w1formal}
Let $W_*(H^m(M,\C))$ be the Deligne weight filtration \cite{D} 
on the cohomology with complex coefficients of a connected 
smooth quasi-projective variety $M$. It follows from a basic 
result of J.~Morgan \cite[Corollary 10.3]{M} that 
$\pi _1(M)$ is $1$-formal if $W_1(H^1(M,\C))=0$. 
By \cite[Corollary 3.2.17]{D}, this vanishing property 
holds whenever $M$ admits a non-singular compactification 
with trivial first Betti number. As noted in \cite{K}, these 
two facts together establish the $1$-formality of fundamental 
groups of complements of projective hypersurfaces. 
\end{example}

\begin{example} 
\label{ex:w3formal}
If $b_1 (G)=0$, it follows from \cite{Q} that 
$E_G \cong \wh (G)=0$, therefore $G$ is $1$-formal. 
If $G$ is finite, Serre showed in \cite{Se} that $G$ is 
the fundamental group of a smooth complex projective variety.
\end{example}

\subsection{On the Malcev completion of a $1$-formal group}
\label{subsec=1malc}

The following test for $1$-formality is well-known, and will be useful 
in the sequel.  For the benefit of the reader, we include a proof.

\begin{lemma} 
\label{lem=quadr}
A finitely presented group $G$ is $1$-formal if and only if the 
Malcev Lie algebra $E_G$ is isomorphic, as a filtered Lie algebra, 
to the completion with respect to degree of a quadratic Lie algebra.
\end{lemma}

\begin{proof} 
Write  $E_G = {\widehat C}$, with $C^*$ a graded Lie algebra presented 
as  $C^*=\bL ^*(X)/\langle Y\rangle $, with $Y \subset \bL ^2(X)$ a 
$\K$-vector subspace. It is enough to prove that $C\cong \h (G)$, 
as graded Lie algebras.   

From \cite{S}, we know that that the image of  
$\partial _G= \partial _{K(G,1)} $ equals the kernel of the Lie bracket 
map $\bigwedge ^2 \gr^1(G) \otimes \K  \to \gr^2(G) \otimes \K$.  
From \eqref{eq=griso}, we get a graded Lie algebra isomorphism 
$\gr^*(G) \otimes \K \cong C^*$.
Via this isomorphism, $H_1(G,\K)=\gr^1(G)\otimes \K$ is identified to 
$X$ and $\im (\partial _G)$ to $Y$. Thus, $C^*\cong \h ^* (G)$. 
\end{proof}

It is convenient to normalize the Malcev completion of a 
$1$-formal group in the following way.  Recall that, for a 
group $G$ with holonomy Lie algebra $\h(G)$, we 
have $\gr_F^*( \widehat { \h }(G))=\h ^*(G)$.  

\begin{lemma} 
\label{lem=norm}
If $G$ is a $1$-formal group, then there is a Malcev completion 
homomorphism
\[
\kk \colon G \to \eexp( \wh (G)) ,
\]
inducing the identity on $\gr^1(G) \otimes \K=\gr^1_F(\wh (G))=\h^1 (G)$.
\end{lemma}

\begin{proof}
Set $\h =  \h (G)$,  $X=G_{\ab}  \otimes \K$ and $Y=\im(\partial_G)$. 
Let $\kk_G \colon G \to \eexp( \widehat { \h })$ 
be a Malcev completion, and set $\phi= \gr^1(\kk_G)\colon X \to X$.  
It follows from \eqref{eq=griso} that $\phi$  is an isomorphism.  
Moreover, $\bL^2( \phi)\colon\bigwedge ^2X \to \bigwedge ^2X$ 
is also an isomorphism, sending $K_1:=\ker ([ \: , \: ]\colon 
\bigwedge ^2X \to \gr^2(G)  \otimes \K)$ onto 
$K_2:=\ker ([ \: , \: ]\colon\bigwedge ^2X \to \h ^2)$. 

As noted in the proof of Lemma \ref{lem=quadr}, $K_1=Y$; 
on the other hand, $Y=K_2$, by  \eqref{eq=holm}. 
Hence, $\phi$ extends to a graded Lie algebra automorphism 
$\phi\colon \h \to \h$.  Upon completion, we obtain a group automorphism 
$\widehat { \phi}\colon \eexp( \widehat { \h }) \to  \eexp( \widehat {\h})$.
The desired normalized Malcev completion is 
$\kk =\widehat { \phi}^{-1} \circ \kk_G$.
\end{proof}

\section{Alexander invariants of $1$-formal groups} 
\label{sec=alex}

Our goal in this section is to derive a relation between the 
Alexander invariant and the holonomy Lie algebra of a 
finitely presented, $1$-formal group.

\subsection{Alexander invariants} 
\label{subsec=alex} 
 
Let $G$ be a group.  Consider the exact sequence
\begin{equation} 
\label{eq=gsq}
\xymatrix{0 \ar[r]& G'/G''  \ar[r]^{j}& 
G/G''   \ar[r]^{p}&   G_{\ab}  \ar[r]&  0},
\end{equation}
where $G'=(G,G)$, $G''=(G',G')$ and $ G_{\ab}=G/G'$. 
Conjugation in $ G/G''$ naturally makes $G'/G''$ into a 
module over the group ring $\Z G_{\ab}$.   We call this module, 
\[
B_G= G'/G'',
\]
the {\em Alexander invariant}\/ of $G$.  If $G=\pi_1(M)$, where $M$ 
is a connected CW-complex, one has the following useful topological 
interpretation for the Alexander invariant. Let $M' \to M$ be 
the Galois cover corresponding to $G' \subset G$. Then
$B_G  \otimes \K=H_1(M',\K)$, and the action of $G_{\ab}$ 
corresponds to the action in homology of the group of
covering transformations.  See \cite{Fox, Cro, Ma} 
as classical references, and \cite{CS1, DM} for 
a detailed treatment in the case of complements of 
hyperplane arrangements and hypersurfaces, respectively.

Now assume the group $G$ is finitely presented. Then 
$B_G  \otimes \K$ is a finitely generated module over the 
Noetherian ring $\K G_{\ab}$. Denote by $I \subset \K G_{\ab}$ 
the augmentation ideal and set $X:= G_{\ab} \otimes \K$. 
The $I$-adic completion $\widehat { B_G  \otimes \K}$ is 
a finitely generated module over 
$\widehat {\K G_{\ab}}=\K [[X]]$, the formal power series ring 
on $X$. Note that $\K [[X]]$ is also the $(X)$-adic completion 
of the polynomial ring $\K [X]$.

Another invariant associated to a finitely presented group $G$ 
is the {\em infinitesimal Alexander invariant}, $B_{\h (G)}$.   
By definition, this is the finitely generated module over $\K[X]$ 
with presentation matrix
\begin{equation} 
\label{eq=infalex}
\nabla:= \delta _3 + \id \otimes \partial _G\colon 
\K [X] \otimes\Big(\bigwedge\nolimits^3X \oplus Y\Big) \to 
\K [X] \otimes \bigwedge\nolimits^2X ,
\end{equation}
where $Y=H_2(G,\K)$ and $ \delta _3 (x \wedge y  \wedge z)
= x\otimes y   \wedge z -y \otimes x  \wedge z +z 
\otimes x \wedge y$; see \cite[Theorem 6.2]{PS2}, and also 
\cite{CS2, MS1} for related definitions.  As the 
notation indicates, the module $B_{\h (G)}$ depends only on 
the holonomy Lie algebra of $G$.

\subsection{The exponential base change} 
\label{subsec=expchange} 

Set $\h = \h (G)$  and denote by  $\h '$ (resp.  $\h ''$) the first 
(resp. the second) derived Lie subalgebras. Consider the exact sequence 
of graded Lie algebras
\begin{equation} 
\label{eq=hsq}
\xymatrix{ 0 \ar[r]&  \h ' / \h '' \ar[r]^{\iota} & \h/\h '' 
\ar^{\pi}[r] & \h _{\ab}\ar[r]&  0 } ,
\end{equation}
where $ \h _{\ab} = \h  / \h' =X.$ Note that the adjoint representation of 
$ \h  / \h ''$ induces a natural graded module structure on $\h'/\h''$, 
over the universal enveloping algebra $U(\h _{\ab})=\K[X]$.

Let $\kk \colon G \to \eexp(\widehat {\h})$ be a normalized 
Malcev completion as in Lemma \ref{lem=norm}. 
Theorem 3.5 from \cite{PS2}  guarantees the factorization 
of $\kk$ through Malcev completion homomorphisms
$\kk '' \colon G/G'' \to \eexp( \widehat { \h/ \h ''})$ and 
$\kk ' \colon G _{\ab}  \to \eexp( \widehat { \h _{\ab}   })$. 
From the definitions, $\kk' \circ p= \widehat {\pi} \circ \kk ''$. 
Furthermore, the completion of the exact 
sequence \eqref{eq=hsq} with respect to degree 
filtrations is still exact. Therefore, there is an induced 
group morphism 
$\kk _0 \colon  G'/G'' \to \eexp( \widehat {\h' / \h ''})$ 
that fits into the following commuting diagram
\begin{equation}
\label{eq=k0diagram}
\xymatrix{
0\ar[r]& G'/G'' \ar^{j}[r] \ar^{\kk_0}[d] &  G/G''\ar^{\kk''}[d]  
\ar[r]^{p}&   G_{\ab}\ar^{\kk'}[d] \ar[r]&0 \\
0\ar[r]&  \eexp( \widehat {\h' / \h ''})  \ar^{\widehat {\iota} }[r] 
&\eexp( \widehat {\h / \h ''}) 
\ar^{\widehat {\pi} }[r] &\eexp( \widehat {\h_{\ab}}) \ar[r]&0
}
\end{equation}

Note that $\h^*_{\ab}=\h ^1_{\ab}=X$ has the structure of an 
abelian Lie algebra. Therefore, by the Campbell-Hausdorff 
formula \eqref{eq=ch}, the exponential group structure on 
$\eexp( \widehat { \h _{\ab} })=X$ coincides with the 
underlying abelian group structure of the $\K$-vector 
space $X$. Moreover, by the normalization property, 
$\kk '$ coincides with the canonical map $G_{\ab} \to X$, 
$a \mapsto a \otimes 1$. Note also that the Lie algebra 
$\h'/\h''$ is abelian, 
hence $ \eexp( \widehat {\h'/\h''})$ is the underlying abelian 
group structure of the $\K$-vector space $\widehat {\h'/\h ''}$, 
again by \eqref{eq=ch}.

Now consider the morphism of $\K$-algebras
\begin{equation} 
\label{eq=e}
\exp\colon \K G_{\ab} \to \K [[X]] , 
\end{equation}
given by $\exp (a)=e^{a \otimes 1}$, for $a \in G_{\ab}$.  
By the universality property of completions, 
there is an isomorphism of filtered  $\K$-algebras, 
$\widehat {\exp}\colon \widehat {\K G_{\ab}} 
\xrightarrow{\,\simeq\,} \K [[X]]$, such that 
the morphism $\exp$ decomposes as 
$\K G_{\ab} \to \widehat{\K G_{\ab}} 
\xrightarrow{\widehat{\exp}} \K [[X]]$, 
where the first arrow is the natural map to the $I$-adic completion.
Note that the completion $\widehat {\h'/\h''}$ has a canonical
module structure over the $(X)$-adic 
completion $\widehat{\K [X]}=\K [[X]]$.

\begin{lemma}
\label{lem=eq}
The map $\kk _0 \otimes \K \colon (G'/G'')\otimes \K \to 
\widehat {\h'/\h''}$ is $\exp$-linear, that is,
\[
(\kk _0 \otimes \K) (\alpha\cdot \beta)= \exp (\alpha)\cdot
(\kk _0 \otimes \K) (\beta)\, ,
\]
for $\alpha \in \K G_{\ab}$ and $\beta\in (G'/G'')\otimes \K$.
\end{lemma}

\begin{proof}  
It is enough to show that 
$\kk ''(a \cdot j(b)  \cdot a ^{-1})=e^{p(a) \otimes 1} \cdot \kk ''(j(b))$ 
for $a \in G/G''$ and $b \in G'/G''$. To check this equality, recall the 
well-known conjugation formula in exponential groups 
(see Lazard \cite{L}), which in our situation says
\[
xyx ^{-1}=\exp(\ad_x)(y) ,
\]
for $x,y \in \eexp( \widehat {\h / \h ''})$. 
Hence $\kk ''(a \cdot j(b)  \cdot a ^{-1})=e^{\ad_{\kk ''(a)}}(\kk ''(j(b)))$, 
which equals $e^{p(a) \otimes 1} \cdot \kk ''(j(b))$, since
$\widehat {\pi}\circ \kk'' (a)=\kk ' \circ p(a)= p(a) \otimes 1$.
\end{proof}

\subsection{Completion of the Alexander invariant} 
\label{subsec=compalex} 

Recall that  $B_G \otimes \K$ is a module over $\K G_{\ab}$, 
and $ \widehat { \h' / \h ''   }$ is a module over 
$\widehat { \K[X]}=\K[[X]]$. 
Shift the canonical degree filtration on $\widehat {\h' / \h ''}$ 
by setting $F'_q\,  \widehat {\h' / \h ''} :=F_{q+2}\, \widehat {\h' / \h ''}$, 
for each $q\ge 0$. 

\begin{lemma} 
\label{lem=filt}
The $\K$-linear map
$\kk _0 \otimes \K\colon B_G \otimes \K \to \widehat {\h' / \h ''}$ 
induces a filtered $\widehat{\exp}$-linear isomorphism between 
$\widehat { B_G \otimes \K}$, endowed with the 
filtration coming from the $I$-adic completion, and 
$\widehat{\h' / \h ''}$, endowed with the shifted 
degree filtration $F'$.
\end{lemma}

\begin{proof} 
We start by proving that
\begin{equation} 
\label{eq=kkinc}
\kk _0  \otimes \K\,(I^q B_G \otimes \K) \subset F'_q \, \widehat {\h'/\h ''}
\end{equation}
for all $q \geq 0$. First note that $ \widehat { \h' / \h ''}= 
F_2  \widehat {\h'/\h''}$, since $\h'$ consists of elements of 
degree at least $2$.  Next, recall that $\exp(I) \subset (X)$. 
Finally, note that
$(X)^rF_s\, \widehat {\h'/\h''} \subset F_{r+s}\, \widehat{\h'/\h''}$
for all $r,s \geq 0$.   These observations, together with 
the $\exp$-equivariance property from Lemma \ref{lem=eq}, 
establish the claim.

In view of  \eqref{eq=kkinc}, $\kk _0  \otimes \K$ induces 
a filtered, $\widehat{\exp}$-linear  map  
from $\widehat { B_G \otimes \K}$ to $\widehat { \h' / \h ''}$.  
We are left with checking  this map is a filtered isomorphism. 
For that, it is enough to show
\begin{equation} 
\label{eq=g1}
\gr^q(\kk _0  \otimes \K)\colon \gr^q_I( B_G \otimes \K) \to  
( \h' / \h '')^{q+2}
\end{equation}
is an isomorphism, for each $q \geq 0$. 

Recall from \eqref{eq=k0diagram} that 
$\kk '' \circ j= \widehat {\iota} \circ \kk _0$. 
By a result of W.~Massey \cite[pp.~400--401]{Ma}, the 
map $j\colon G'/G''\to G/G''$ induces isomorphisms
\begin{equation} 
\label{eq=g2}
\gr^q(j)  \otimes \K \colon \gr^q_I( B_G \otimes \K) 
\xrightarrow{\,\simeq\,}   \gr^{q+2}  (G / G '') \otimes \K  
\end{equation}
for all $q \geq 0$. Since $\kk''$ is a 
Malcev completion, it induces isomorphisms
\begin{equation} 
\label{eq=g3}
\gr^q(\kk '')\colon  \gr^{q  }  (G / G '') \otimes 
\K \xrightarrow{\,\simeq\,} ( \h / \h '')^{q }   
\end{equation}
for all $q \geq 1$.  Finally (and evidently), the inclusion map 
$\iota\colon  \h' / \h '' \to  \h / \h ''$ induces isomorphisms 
$\gr^q(\iota)$, and thus $\gr^q (\widehat {\iota})$, for all $q\ge 2$.  
This finishes the proof.
\end{proof}

\begin{lemma} 
\label{lem=end}
The $\K[[X]]$-module $\widehat {B_{\h (G)}}$ is filtered 
isomorphic to the module $\widehat { \h' / \h ''}$, 
endowed with the shifted filtration $F'$.
\end{lemma}

\begin{proof}  
Assign degree $q$ to $\bigwedge ^qX$ and degree $2$ to 
$Y$ in \eqref{eq=infalex}. Then $ \h' / \h ''$, viewed as a 
graded  $\K[X]$-module via the exact sequence \eqref{eq=hsq}, 
is graded isomorphic to the  $\K[X]$-module $\coker (\nabla)$, 
see  \cite[Theorem 6.2]{PS2}. Taking $(X)$-adic completions, 
the claim follows.
\end{proof}

Putting Lemmas  \ref{lem=filt} and \ref{lem=end} together, 
we obtain the main result of this section.

\begin{theorem} 
\label{thm=complalex}
Let $G$ be a $1$-formal group. Then the $I$-adic completion of 
the Alexander invariant, $\widehat { B_G  \otimes \K}$, 
is isomorphic to the $(X)$-adic  completion of the 
infinitesimal Alexander invariant, $\widehat { B_{\h (G)}}$,
by a filtered $\widehat {\exp}$-linear isomorphism.
\end{theorem}

\subsection{Fitting ideals}
\label{subsec=fitt}

We now recall some basic material from \cite[\S20.2]{E}. 
Let $R$ be a commutative Noetherian ring, and let $N$ be 
a finitely generated $R$-module.  Then $N$ admits a presentation 
of the form 
\begin{equation} 
\label{eq=npres}
\xymatrix { R^r \ar^{\nabla }[r] & R^s \ar[r] & N \ar[r] &  0  } .   
\end{equation}
Define the $k$-th Fitting ideal of $N$ as follows:  $\F_k(N)$ is 
the ideal generated by the $(s-k)$-minors 
of $\nabla$, if $0<s-k \leq \min\{r,s\}$,  $\F_k(N)=R$ if 
$s-k \leq 0$, and $\F_k(N)=0$, otherwise.
The Fitting ideals form an ascending chain, independent of 
the choice of presentation for $N$. Their construction commutes 
with base change. 

If $R=\C[A]$ is the coordinate ring of an affine variety, 
we may consider the decreasing filtration of $A$ by the 
{\it Fitting loci} of $N$, 
\begin{equation} 
\label{eq=floci}
\W_k(N):=Z(\F_{k-1}(N)) \subset A .
\end{equation}
These loci (which depend only on the coherent 
sheaf $\widetilde N$) define a decreasing filtration of $A$ by 
closed reduced subvarieties.

The next (elementary) lemma will be useful in the sequel.

\begin{lemma} 
\label{lem=useful}
Let $A$ be an affine variety, with coordinate ring $R$. For a point 
$t\in A$, denote by $\m_t \subset R$ the corresponding maximal 
ideal.  Let $N$ be a finitely generated  $R$-module. Then
\[
t \in \W_k(N) \iff \dim _{\C}(N/\m_tN) \geq k .
\]
\end{lemma}

\begin{proof} 
First note that $N/\m_tN=R/\m_t\otimes _R N$. Assuming 
$0<s-k+1 \leq \min\{r,s\}$ in presentation  \eqref{eq=npres}, 
we may translate definition \eqref{eq=floci} to
\begin{equation} 
\label{eq=n1}
t \in \W_k(N) \iff \rank _{\C}(R/\m_t \otimes_R \nabla) \leq s-k . 
\end{equation}
We also have 
\begin{equation} 
\label{eq=n2}
\rank _{\C}(R/\m_t\otimes_R \nabla) \leq s-k \iff  
\dim _{\C}(R/\m_t \otimes _R N) \geq k  . 
\end{equation}
These two equivalences yield the claim in this case. If $s-k+1>r$, 
then $\W_k(N)=A$ by definition.  In this situation, we also have 
$\dim _{\C}(R/\m_t\otimes _R N) \geq s-r \geq k.$ 
The remaining cases are similar.
\end{proof} 

\subsection{Fitting loci and $1$-formality}
\label{subsec=fitt1f}

Let $G$ be a finitely presented group, with Alexander invariant 
$B_G$ and infinitesimal Alexander invariant $B_{\h(G)}$. 
Use  the isomorphism $\widehat{\exp} \colon \widehat{\C G_{\ab}}
\xrightarrow{\simeq} \C[[X]]$ to identify $\widehat {\C G_{\ab}}$, 
the $I$-adic completion of the group algebra $\C G_{\ab}$,  with 
$\C[[X]]$, the $(X)$-adic completion of the polynomial ring $\C[X]$. 

\begin{lemma} 
\label{lem=fext} 
If $G$ is $1$-formal, then
\[
\widehat {\exp}\colon \F_k(B_G\otimes \C)\, \C[[X]] 
\xrightarrow{\,\simeq\,} \F_k(B_{\h(G)})\, \C[[X]] ,
\]
for all $k \geq 0$.
\end{lemma}

\begin{proof} 
This is a direct consequence of Theorem \ref{thm=complalex}, 
by base change.
\end{proof}

Let $T_1 \T_G=\Hom(G,\C)$ be the tangent space at the 
origin $1$ to the complex analytic torus $\T_G=\Hom(G,\C^*)$.  
Then  $T_1\T_G$ is the Lie algebra of $\T_G$, 
and the exponential map $\C\to \C^*$, $a\mapsto e^a$, 
induces a local analytic isomorphism
\begin{equation}
\label{eq=expo}
\exp\colon (T_1\T_G,0) \xrightarrow{\,\simeq\,} (\T_G, 1) .
\end{equation}

\begin{prop} 
\label{prop=floc} 
If $G$ is $1$-formal, then the exponential map induces an 
isomorphism of analytic germs
\[
\exp\colon (\W_k (B_{\h(G)}),0) \xrightarrow{\,\simeq\,}
(\W_k(B_G\otimes \C),1) ,
\]
for all $k \geq 0$.
\end{prop}

\begin{proof}
In geometric terms, the $I$-adic completion map  can be written 
as the composite $\C G_{\ab} \to \C \{X\} \to \C[[X]]$, 
where the first arrow takes global regular functions on $\T_G$ 
to their analytic germs at $1$, and the second arrow is Taylor 
expansion. One has a similar decomposition for the $(X)$-adic 
completion, $\C [X] \to \C \{X\} \to \C[[X]]$, obtained by 
taking germs (respectively, Taylor expansions) at $0$. 

With these identifications, the map $\exp \colon \C G_{\ab} \to 
\C[[X]]$ from \eqref{eq=e} can be viewed as the composite 
$\C G_{\ab} \to \C \{X\} \to \C \{X\} \to \C[[X]]$, where the 
middle arrow is the map induced by \eqref{eq=expo} on 
local coordinate rings. 

Since $\C[[X]]$ is faithfully flat over $\C \{X\}$, see \cite[p.~36]{T}, 
we obtain from Lemma \ref{lem=fext} an exponential identification
\[
\exp\colon \F_k(B_G\otimes \C)\, \C \{X\} \xrightarrow{\,\simeq\,}
\F_k(B_{\h(G)})\, \C \{X\}
\]
for all $k \geq 0$. The claim follows by taking the corresponding 
germs of zero sets.
\end{proof}

\section{Germs of cohomology support loci} 
\label{sec=germs}

In this section, we study the cohomology support loci of a 
finitely presented group $G$: the resonance varieties 
$\R_k(G)$ and the characteristic varieties $\V_k(G)$.  
Under a $1$-formality assumption on $G$, we 
show that the exponential map restricts to an isomorphism 
of analytic germs between $(\R_k(G),0)$ and $(\V_k(G),1)$. 
We work over the field $\K=\C$, unless otherwise specified.

\subsection{Resonance varieties} 
\label{subsec=resvar}

Let $H^1$ and  $H^2$ be two finite-dimensional complex 
vector spaces and let $\mu\colon H^1 \wedge  H^1 \to  H^2$ 
be a linear map. Set $H^0=\C$. For an element $z \in H^1$, 
denote by $\cdot z$ right-multiplication by $z$ in the exterior 
algebra $\bigwedge^{*} H^1$, and by $\mu_z$ the composite 
$H^1\xrightarrow{\cdot z} H^1\wedge H^1 \xrightarrow{\mu} 
H^2 $.  We then have a cochain complex
\begin{equation} 
\label{eq=aom}
(H^{\bullet}, \mu _z)\colon 
\xymatrix{0 \ar[r]& H^0  \ar[r]^{\cdot  z} & H^1    
\ar[r]^{\mu_z } & H^2  \ar[r]&  0}, 
\end{equation}
with cohomology denoted by $H^*(H^{\bullet}, \mu _z)$.
For each $k \geq 0$, the $k$-th {\it resonance variety} of 
$\mu$ is defined by 
\begin{equation} 
\label{eq=resmu}
\R_k(\mu):=\{ z \in H^1 \mid  \dim H^1(H^{\bullet}, \mu _z)   \geq k \}\, .  
\end{equation}
Clearly, the jumping loci  $\{\R_k(\mu)\}_k$ form a decreasing 
filtration of $H^1$ by closed homogeneous subvarieties, starting 
from $\R_0(\mu)=H^1$. 

For a CW-complex $M$ as above, take $H^i=H^i(M,\C)$, 
take $\mu $ to be the cup-product map $\cup_M$, 
and define the corresponding resonance varieties
$\R_k(M):= \R_k( \cup _M  )$. 
Similarly, for a finitely presented group $G$, define 
\begin{equation} \label{eq=resg}
\R_k(G):= \R_k( K(G,1)  )  .  
\end{equation}

Obviously, the resonance varieties \eqref{eq=resmu} depend 
only on the corestriction of $\mu$ to its image. By applying 
this remark to the classifying map $f\colon M \to K(\pi_1(M),1)$, 
which induces an isomorphism on $H^1$ and a monomorphism 
on $H^2$, we see that
\begin{equation} 
\label{eq=reseq}
\R_k(M)= \R_k( \pi_1(M)  )   
\end{equation}
for all $k \geq 0$.

\begin{lemma} 
\label{lem=step1}
Let $G$ be a finitely presented group. Then the equality
\[
\W_k (B_{\h(G)}) \setminus \{0\}=\R_k(G)\setminus \{0\}
\]
holds for all $k \geq 0$.
\end{lemma}

\begin{proof} 
Pick any $z \in H^1(G, \C)\setminus \{0\}$. By applying 
Lemma  \ref{lem=useful} to $B_{\h(G)}$, we infer that 
$z \in \W_k (B_{\h(G)})$ if and only if $
\dim _{\C}\coker (\nabla (z)) \geq k$.  Note that 
$H^1(G, \C)$ equals the dual vector space 
${}^{\sharp}X$ of $X$. Consider the cochain complex
$(\bigwedge ^{\bullet} {}^{\sharp}X, \lambda _z)$, where 
$\lambda _z$ denotes left multiplication by $z$.  Since 
$z \ne 0$, this complex is obviously exact. Let 
$(\bigwedge^{\bullet} X, {}^{\sharp} \lambda _z)$ be the 
dual chain complex, which is again exact. It is straightforward 
to check that the restriction of $ {}^{\sharp} \lambda _z$ to 
$\bigwedge ^{3} X$ equals $\delta _3(z)$, with $\delta _3$ 
as in the presentation  \eqref{eq=infalex} of $B_{\h(G)}$.

Denoting by $\delta _2(z)$ the restriction of $  ^{\sharp} \lambda _z$ 
to $\bigwedge ^2 X$, use the exactness of the complex 
$(\bigwedge ^{\bullet} X,{}^{\sharp} \lambda _z)$ 
to obtain the following isomorphism:
\begin{equation} 
\label{eq=ex}
\coker (\nabla (z)) \cong \im ( \delta _2(z))/ 
\im ( \delta _2(z)\circ \partial _G) .
\end{equation}
By exactness again,  $ \dim _{\C}\im ( \delta _2(z))=n-1$. Hence 
$z \in \W_k (B_{\h(G)})$ if and only if 
$\rank ( \delta _2(z)\circ \partial _G) \leq n-1-k$.  
Using \eqref{eq=holm}--\eqref{eq=holg} and 
\eqref{eq=aom}--\eqref{eq=resg}, we see   
that the linear map dual to 
$\delta _2(z)\circ \partial _G$ is $- \mu _z$. Consequently, 
$z \in \W_k (B_{\h(G)})$ if and only if
$\rank ( \mu_z) \leq n-1-k$, that is, if and only if 
$z \in \R_k(G)$, see \eqref{eq=resmu}--\eqref{eq=resg}.
\end{proof} 

\subsection{Characteristic varieties} 
\label{subsec=charvar}

Let $M$ be a connected CW-complex with finite $2$-skeleton. 
Set $G=\pi_1(M)$ and consider the character torus 
$\T_G=\Hom(G,\C^*)$.  This is an algebraic group of the form 
$(\C^*)^n \times F$, where $n=b_1(M)$ and $F$ is a finite abelian 
group. Identifying the point $\rho \in \T_G$ with a rank one local 
system ${}_{\rho} \C$ on $M$ (that is, a left one-dimensional 
$\C$--representation of $\pi_1(M)$), one may define the $k$-th 
{\em characteristic variety}\/ for all $k \geq 0$ by
\begin{equation} 
\label{eq=charm}
\V_k(M):=\{ \rho \in \T_G \mid  \dim H^1(M, {} _{\rho}\C) \geq k \} .  
\end{equation}
Here $H^{\bullet}(M, {} _{\rho}\C)$ denotes twisted cohomology, 
see e.g.~\cite[Chapter VI]{W}.  Viewing $\rho$ as a right representation, 
it follows by duality that we may replace $H^{1}(M, {} _{\rho}\C)$ 
by twisted homology, $H_{1}(M, {} _{\rho}\C)$, in definition 
\eqref{eq=charm}.

Plainly, $\{\V_k(M)\}_k$ is a 
decreasing filtration of the torus $\T_G$ by closed algebraic 
subvarieties, starting from $\V_0(M)=\T_G$. Describing this 
filtration is equivalent to understanding the degree one 
(co)homology of $M$ with coefficients in an arbitrary rank 
one local system, a very difficult task in general. 

For a finitely presented group $G$, define
\begin{equation} 
\label{eq=charg}
\V_k(G):=\V_k(K(G,1)) .  
\end{equation}
If $M$ is a $CW$-complex as above, with $G=\pi_1(M)$, then 
an Eilenberg-MacLane space $K(G,1)$ can be obtained from 
$M$ by attaching cells of dimension $\ge 3$; hence 
\begin{equation} 
\label{eq=chareq}
\V_k(M)=\V_k(\pi_1(M)) ,  
\end{equation}
for all $k \geq 0$.

Let  $G=\langle x_1,\dots ,x_s \mid w_1,\dots ,w_r \rangle$ 
be a finite presentation for $G$. The $2$-complex $M$ 
associated to this presentation has $s$ one-cells, corresponding 
to the generators $x_i$, and $r$ two-cells, attached according 
to the defining relations $w_j$. Denote by 
\begin{equation} 
\label{eq=eqchains}
{\widetilde C}_{\bullet}= C_{\bullet}(\widetilde{M}) \colon 
\xymatrix{ 0 \ar[r]& (\Z G )^r  \ar^{d_2 }[r]& (\Z G )^s  
\ar^{d_1 }[r]& \Z G \ar^{\epsilon}[r]& \Z }
\end{equation}
the augmented $G$-equivariant cellular chain complex 
of the universal cover ${\widetilde M}$.

\begin{remark} 
\label{rk=semi} 
If we tensor the chain complex \eqref{eq=eqchains} (where the last 
term $\Z$ is replaced by $0$) with $\C$, via the ring extension 
$\Z G \to  \C G_{\ab}  \xrightarrow{\rho} \C$ 
associated to a character $\rho$, we obtain a complex of 
finite-dimensional $\C$-vector spaces, namely
\[
\xymatrix{ 0  \ar[r]& \C ^r   \ar[r]^{d_2(\rho)} & \C ^s  
 \ar[r]^{d_1(\rho)} & \C  \ar[r]& 0 }. 
\]
Here the differentials $d_k(\rho)$ are represented by matrices 
whose entries are regular functions in $\rho \in \T_G$. It follows 
that the function
\[
\dim H_1(G,{}_{\rho}\C)=\dim \ker d_1(\rho) - \rank d_2(\rho)= 
s-\rank d_1(\rho)- \rank d_2(\rho)
\]
is lower semi-continuous with respect to $\rho$, i.e.,
\[
\dim H_1(G,{}_{\rho}\C) \leq \dim H_1(G,{} _{\rho _0}\C)
\]
for $\rho$ in a neighborhood of a fixed character $\rho_0$.
\end{remark}

Let $\Z G \to \C G_{\ab}$ be the base change associated 
to the abelianization morphism.  As mentioned in \S\ref{sec=alex}, 
there is  a natural $\C G_{\ab}$-module identification
\begin{equation} 
\label{eq=cellalex}
B_G \otimes \C = H_1( \C G_{\ab}  \otimes_{\Z G}  {\widetilde C}_{\bullet}) .
\end{equation}

Let $\Z G \to  \C G_{\ab}  \xrightarrow{\rho} \C$ be 
the change of rings corresponding to a character $\rho \in \T_G$. 
As follows from \eqref{eq=chareq},
\begin{equation} 
\label{eq=cellrho}
\rho \in \V_k(G) \iff \dim_{\C}H_1({}_{\rho}\C 
\otimes _{ \Z G} {\widetilde C}_{\bullet} ) \geq k .
\end{equation}
By applying Lemma  \ref{lem=useful} to $N=B_G \otimes \C$ 
and using \eqref{eq=cellalex}, we infer that
\begin{equation} \label{eq=cellw}
\rho \in \W_k(B_G \otimes \C) \iff \dim_{\C}({}_{\rho}\C 
\otimes _{ \C G_{\ab}} H_1K_{\bullet}) \ge k  ,
\end{equation}
where $K_{\bullet}$ denotes the $\C G_{\ab}$-chain complex 
$\C G_{\ab} \otimes_{\Z G} {\widetilde C}_{\bullet}$.

With these preliminaries, we may state the following analogue of 
Lemma \ref{lem=step1}.  

\begin{lemma} 
\label{lem=step2}
Let $G$ be a finitely presented group. Then the equality
\[
\W_k(B_G\otimes \C)\setminus \{1\} = \V_k(G)\setminus \{1\}
\]
holds for all $k\ge 0$.
\end{lemma}

\begin{proof}
Let $R= \C G_{\ab}$. 
By \eqref{eq=cellrho} and \eqref{eq=cellw}, 
it is enough to show 
$H_1({}_{\rho}\C \otimes _R K_{\bullet})= 
{} _{\rho}\C \otimes _R H_1K_{\bullet}$, 
for $\rho \in \T_G\setminus \{1\}$. 
An analysis of the spectral sequence associated to 
the free, finite chain complex $K_{\bullet}$ over the 
ring $R$, and the base change $\rho\colon R \to \C$, 
\[
E^2_{s,t}=\tor^R_s({}_{\rho}\C,H_tK_{\bullet}) \Rightarrow 
H_{s+t}({}_{\rho}\C \otimes _RK_{\bullet}) ,
\]
see \cite[Theorem XII.12.1]{ML}, shows we only need to check  
$\tor^R_*({}_{\rho}\C, {} _{\epsilon}\C)=0 $, where 
$\epsilon\colon R\to \C$ is the augmentation map. 

Denote by $F \subset G_{\ab}$ the torsion subgroup, and by 
$G_{\ab}/F =\Z^n$ the torsion-free part of $G_{\ab}$. Since 
$R \cong \C \Z^n \otimes \C F$, it follows from the K\"{u}nneth 
formula that
\begin{equation*} 
\label{eq=ku}
\tor^R_*({}_{\rho}\C, {} _{\epsilon}\C)=\tor^{ \C \Z^n }_*({}_{\rho'}\C, 
{}_{\epsilon}\C) \otimes
\tor^{ \C F }   _*({}_{\rho''}\C, {} _{\epsilon}\C) , 
\end{equation*}
where $\rho=(\rho',\rho'')\in \T_{\Z^n}\times \T_F$. It is well-known that
$\tor^{ \C \Z^n }_*({}_{\rho'}\C, {} _{\epsilon}\C)=0$, if $\rho'\ne 1$, 
see e.g.~\cite[pp.~215--216]{Ha}. Thus, it is enough to show 
$\tor^{ \C F}_*({}_{\rho''}\C, {} _{\epsilon}\C)=0$,  if $\rho''\ne 1$, 
for any finite cyclic group $F$.  This can be easily checked, using 
the standard periodic resolution of ${}_{\epsilon}\C$ over $\C F$, 
see e.g.~\cite[IV.7]{ML}.
\end{proof}

\subsection{Matching analytic germs} 
\label{subsec=expmap}

The next theorem proves the first part of 
Theorem \ref{thm=tcfintro} from the Introduction.  
It implies, in particular, that the twisted (degree one) 
cohomology with rank one local systems is computable, 
near the trivial character $1$ and under a formality 
assumption, only in terms of the cup-product map 
on $H^1$.

\begin{theorem}
\label{thm=analiso}
Let $G$ be a $1$-formal group. Then the exponential map 
from $(T_1\T_G,0)$ to $(\T_G, 1)$ 
restricts for each $k \geq 0$ to an isomorphism of analytic germs,  
\[
\exp\colon (\R_k(G),0)  \xrightarrow{\,\simeq\,} (\V_k(G),1) . 
\]
\end{theorem}

\begin{proof}
We need to do two things: 
replace in Proposition \ref{prop=floc} the germ 
$(\W_k (B_{\h(G)}),0)$ by $(\R_k(G),0)$, and the 
germ $(\W_k (B_G \otimes \C),1)$ by 
$(\V_k(G),1)$. 

Set $n=b_1(G)$. Then plainly $0 \in \R_k(G)$ and 
$1 \in \V_k(G)$ for $k \leq n$. It is also clear that
$(\R_k(G),0)$ is the  germ of the  empty set  for $k>n$. 
To see that the same holds for the germ $(\V_k(G),1)$, 
note that $\dim _{\C}H^1(G, {}_{\rho}\C) \le 
\dim _{\C}H^1(G, {} _{1}\C)=n$ for any $\rho$ near $1$, 
by semi-continuity, see Remark \ref{rk=semi}. 
Therefore it is enough to make the two aforementioned 
replacements only away from $0$ (resp.~$1$).
This can be done using Lemmas \ref{lem=step1} and 
\ref{lem=step2}.
\end{proof}

\begin{remark}
\label{rem=arr}
Suppose $M$ is the complement of a hyperplane arrangement 
in $\C^{m}$, with fundamental group $G=\pi_1(M)$.  
In this case, Theorem \ref{thm=analiso} can be deduced from 
results of Esnault--Schechtman--Viehweg \cite{ESV} and 
Schechtman--Terao--Varchenko \cite{STV}.  In fact, one 
can show that there is a combinatorially defined open 
neighborhood $U$ of $0$ in $H^1(M, \C)$ with the 
property that  
$H^*(M,{}_{\rho}\C) \cong H^*(H^{\bullet}(M,\C),  \mu_z)$, 
for all $z\in U$,  where $\rho=\exp (z)$. 
A similar approach works as soon as $W_1 (H^1 (M,\C))=0$,  
in particular, for complements of arrangements 
of hypersurfaces in projective or affine space.   
For details, see \cite[Corollary 4.6]{DM}.
\end{remark}

The local statement from Theorem \ref{thm=analiso} is the 
best one can hope for, as shown by the following classical example.

\begin{example} 
\label{ex:knots}
Let $M$ be the complement in $S^3$ of a tame knot $K$. 
Since $M$ is a homology circle, it follows easily that $M$ is 
a formal space; therefore, its fundamental group,  
$G=\pi_1(M)$, is a $1$-formal group. Let $\Delta(t)\in \Z [t, t^{-1}]$ 
be the Alexander polynomial of $K$. It is readily seen that 
$\V_1(G)= \{ 1\} \coprod \Zero (\Delta)$ and $\R_1(G)= \{ 0\}$. 
Thus, if $\Delta(t)\not\equiv 1$, then $\exp(\R_1(G)) \ne \V_1(G)$. 

Even though the germ of $\V_1(G)$ at $1$ provides no information 
in this case, the global structure of $\V_1(G)$ is quite meaningful. 
For example, if $K$ is an algebraic knot, then  $\Delta(t)$ must be  
product of cyclotomic polynomials, as follows from work of 
Brauner and Zariski from the 1920s, see \cite{Mi}.
\end{example}

\begin{remark} 
\label{rem:curves}
Let $M$ be the complement in $\C^2$ of an 
algebraic curve, with fundamental group $G=\pi_1(M)$, 
and let $\Delta(t)\in \Z [t, t^{-1}]$ be the Alexander polynomial 
of the total linking cover, as defined by Libgober; see \cite{Li94} 
for details and references. It was shown in \cite{Li94} that 
all the roots of $\Delta(t)$ are roots of unity. This 
gives restrictions on which finitely presented groups can be 
realized as fundamental groups of plane curve complements. 

Let $\Delta^G \in \Z [t_1^{\pm 1},\dots, t_n^{\pm 1}]$ be the
multivariable Alexander polynomial of an arbitrary quasi-projective 
group.  Starting from Theorem \ref{thm=posobstr}\eqref{a2}, we 
prove in \cite{DPS-codone} that $\Delta^G$ must have a single 
essential variable, if $n\ne 2$. Examples from \cite{P07} show 
that this new obstruction efficiently detects non-quasi-projectivity 
of (local) algebraic link groups.  Note that all roots of the 
one-variable (local) Alexander polynomial $\Delta (t)$ of 
algebraic links are roots of unity; see \cite{Mi}.
\end{remark}

\subsection{Initial ideals and tangent cones} 
\label{sec=init}

Let $(X,0)$ be a reduced 
analytic space germ at the origin of $\C^n$ and let 
$I=I(X,0) \subset \OO _n$ be the ideal of analytic function 
germs at the origin of $\C^n$ vanishing on $(X,0)$. 
Any non-zero $f \in \OO _n$ can be uniquely written as a sum
\[
f=f_m + f_{m+1}+ \cdots ,
\]
where each $f_k$ is a homogeneous polynomial of 
degree $k$, and $f_m \ne 0$.   We call $f_m$ the 
initial form of the germ $f$ and denote it by $\init (f)$. 

The  {\em initial ideal}\/ of $I$, denoted 
$\init (I)$, is the polynomial ideal  spanned by all 
initial forms of non-zero elements of $I$. The {\em tangent cone}\/ 
of the germ $(X,0)$, denoted  $TC_0(X)$, is the affine cone 
in $\C^n$ given by the zero-set of the initial ideal $\init (I)$. If $X$ 
is an algebraic subvariety in $\C^n$ and $p \in X$ is any point,
then the tangent cone $TC_p(X)$ is defined in the obvious way, i.e., 
by translating the germ $(X,p)$ to the origin. For more geometric 
details on this definition see Whitney \cite[pp.~210--228]{Wy};  
for relations to Gr\"{o}bner bases, see Cox-Little-O'Shea  \cite{CLO}. 

Tangent cones enjoy the following functoriality property. 
Let $(X,0) \subset (\C^n,0)$ and  $(Y,0) \subset (\C^p,0)$ 
be two reduced analytic space germs, and let 
$f\colon (\C^n,0) \to (\C^p,0)$ be an analytic map germ 
such that $f(X,0) \subset (Y,0)$. Then the differential
$d_0f\colon T_0\C^n \to T_0\C^p$ satisfies 
$d_0f(TC_0(X)) \subset TC_0(Y)$.  In particular, 
if $f$ is a local analytic isomorphism, then its 
differential induces a linear isomorphism, 
$d_0f\colon TC_0(X) \xrightarrow{\,\simeq\,} TC_0(Y)$.

\subsection{The tangent cone formula} 
\label{sec=tgcone}
  
Using Theorem \ref{thm=analiso} and the above discussion, 
we obtain the following tangent cone formula, which generalizes 
results from \cite{CS2, Li01}, and finishes the proof of 
Theorem \ref{thm=tcfintro} from the Introduction. 

\begin{theorem} 
\label{thm=tcf} 
If the group $G$ is $1$-formal, then $TC_1(\V_k(G))=\R_k(G)$, 
for all $0 \leq k \leq b_1(G)$.
\end{theorem}

\begin{remark}
\label{rem=compu}
In general, it is difficult to verify the $1$-formality of a finitely 
presented group $G$ directly from Definition \ref{def=1formal}. 
The above tangent cone formula provides a new, computable 
obstruction to $1$-formality. Indeed, recall the equivariant 
chain complex ${\widetilde C}_{\bullet}$ from  \eqref{eq=eqchains}. 
It follows from  \eqref{eq=cellrho} that the characteristic 
varieties $\V_k(G)$ may be computed from the Fitting ideals 
of the $\C{G}_{\ab}$-module presented by the {\it Alexander matrix} 
associated to $\C{G}_{\ab} \otimes_{\Z{G}}d_2$.

It is a classical fact that the Alexander matrix is computable 
directly from a finite presentation for $G$, by means of the 
Fox calculus, see \cite{Fox}. The passage from the variety 
$\V_k(G)$ to the tangent cone $TC_1(\V_k(G))$ is achieved by 
effective commutative algebra methods, as described in \cite{CLO}; 
for details, see e.g.~\cite{CS1, CS2}.

It is even simpler to determine the resonance varieties $\R_k(G)$:  
compute the cup-product map $\cup _G$ and the 
holonomy Lie algebra $\h (G)$ directly from the 
group presentation, by Fox calculus (or, equivalently, 
by Magnus expansion), and then use Fitting ideals 
to compute $\R_k(G)$ as in Lemma  \ref{lem=step1}; 
for details, see e.g.~\cite{MS1}. 
\end{remark}

\section{Regular maps onto curves} 
\label{sec=curves}

In this section, we discuss the relationship between 
the cohomology jumping loci of a quasi-compact K\"{a}hler 
manifold $M$, holomorphic maps from $M$ to complex curves, 
and isotropic subspaces in $H^*(M,\C)$. The basic tool is 
a result of Arapura \cite{A}, which we start by recalling.

\subsection{Arapura's theorem}
\label{subsec=arapura}
First, we establish some terminology.   By a {\em curve}\/  
we mean a smooth, connected, complex algebraic variety of 
dimension $1$.  A  curve $C$ admits a canonical compactification 
$\oC$, obtained by adding a finite number of points. 

Following  \cite[p.~590]{A}, we say a map  $f\colon M \to C$ from 
a connected, quasi-compact K\"{a}hler manifold $M$ to a 
curve $C$ is {\em admissible}\/ if $f$ is holomorphic
and surjective, and has a holomorphic, surjective extension with 
connected fibers, $\overline{f}\colon \oM \to \oC$, where $\oM$ 
is a smooth compactification, obtained by adding divisors with 
normal crossings.

With these preliminaries, we can state Arapura's result 
\cite[Proposition V.1.7]{A}, in a slightly modified form, 
suitable for our purposes.

\begin{theorem} 
\label{thm:vadm}
Let $M$ be a connected, quasi-compact K\"{a}hler manifold. 
Denote by $\{ \V^{\alpha}\}_{\alpha}$ the set of irreducible components 
of $\V_1(\pi_1(M))$ containing $1$. If  $\dim \V^{\alpha}>0$, then 
the following hold.

\begin{enumerate}
\item \label{va1} 
There is an admissible map, 
$f_{\alpha}\colon M \to C_{\alpha}$, where $C_{\alpha}$ 
is a smooth curve with $\chi (C_{\alpha})<0$, 
such that  
\[
\V^{\alpha}=f_{\alpha}^* \T_{\pi_1(  C_{\alpha} ) }
\]
and $(f_{\alpha}) _{\sharp}\colon \pi_1(M) \to \pi_1(  C_{\alpha}  )$ 
is surjective.

\item \label{va2}
There is an isomorphism
\[
H^1(M, {}_{f_{\alpha}^* {\rho}}\C)\cong H^1( C_{\alpha}, {}_{\rho}\C) ,
\]
for all except finitely many local systems 
${\rho} \in \T _{\pi_1(  C_{\alpha} )}$.
\end{enumerate}
\end{theorem}

\begin{proof}
Part~\eqref{va1}.  
Proposition V.1.7 from \cite{A} guarantees the existence of an 
admissible map $f_{\alpha}\colon M \to C_{\alpha}$ to a curve 
$C_{\alpha}$ with $\chi (C_{\alpha})<0$  
such that  $\V^{\alpha}=f_{\alpha}^* \T_{\pi_1(C_{\alpha})}$. 
The surjectivity of $(f_{\alpha}) _{\sharp}$ is implicit 
in Arapura's proof.  It follows from two standard facts. First, 
the generic fibers of $f_{\alpha}$ are connected, as soon as this 
happens for $\overline{f}_{\alpha}$, and second, there is a 
Zariski open subset $U_{\alpha}\subset C_{\alpha}$ with the 
property that 
$f_{\alpha} \colon f_{\alpha}^{-1}(U_{\alpha})\to U_{\alpha}$ 
is a locally trivial fibration.

Part~\eqref{va2}. This is stated in the proof of Proposition V.1.7 
from \cite{A}, with ``all except finitely many'' replaced by 
``infinitely many,'' but a careful look at Arapura's argument 
reveals that actually the finer property holds.
\end{proof}

When $M$ is compact, similar results to Arapura's were obtained 
previously by Beauville \cite{Beau2} and Simpson \cite{Sim}.
The closely related construction of regular mappings from 
an algebraic variety $M$ to a curve $C$ starting with suitable 
differential forms on $M$ goes back to Castelnuovo--de Franchis, 
see Catanese \cite{Cat}. When both $M$ and $C$ are compact, 
the existence of a non-constant holomorphic map $M \to C$ 
is closely related to the existence of an epimorphism 
$\pi_1(M) \surj \pi_1(C)$, see Beauville \cite{Beau1}
and Green--Lazarsfeld \cite{GL}. In the non-compact case, 
this phenomenon is discussed in Corollary V.1.9 from \cite{A}.

\begin{corollary} 
\label{cor:radm}
Let $M$ be a connected, quasi-compact K\"{a}hler manifold, 
with fundamental group $G=\pi_1(M)$.  
Assume that $G$ is $1$-formal, with $b_1(G)>0$ and 
$\R_1(G)\ne \{0\}$. Then all irreducible components of 
$\V_1(G)$ containing $1$ are positive-dimensional.
Realize each irreducible component $\V^{\alpha}$ of $\V_1(G)$ 
containing $1$ by pullback via $f_{\alpha}\colon M \to C_{\alpha}$, 
as in Theorem \ref{thm:vadm}\eqref{va1}. Then 
$T_1(\V^{\alpha})= f_{\alpha}^*H^1( C_{\alpha},\C)$,
\begin{equation} 
\label{eq=fir}
\R_1(G)=\bigcup\nolimits _{\alpha}f_{\alpha}^*H^1(C_{\alpha},\C)  ,
\end{equation}
and this decomposition coincides with the decomposition 
of $\R_1(G)$ into irreducible components, 
$\R_1(G)=\bigcup_{\alpha}\R^{\alpha}$.  In particular, 
$\dim \R^{\alpha} =b_1(C_{\alpha})$, for all $\alpha$.
\end{corollary}

\begin{proof}
Since $b_1(G)>0$ and $\R_1(G)\ne \{0\}$,  all irreducible 
components $\R^{\alpha}$ of the algebraic set $\R=\R_1(G)$ 
are positive dimensional. Since $G$ is $1$-formal, the 
tangent cone formula from Theorem \ref{thm=tcf} implies 
that all irreducible components $\V^{\alpha}$ of $\V_1(G)$ 
containing $1$ are positive dimensional. Hence, 
Theorem \ref{thm:vadm}\eqref{va1} applies.  

By Theorem \ref{thm=tcf}, 
$\R_1(G)=TC_1(\V_1(G))=\bigcup\nolimits _{\alpha}TC_1(\V^{\alpha})$.
From Theorem \ref{thm:vadm}\eqref{va1}, we know that each 
$\V^{\alpha}$ is an algebraic connected subtorus of $\T_G$, 
isomorphic to $\T_{\pi_1(  C_{\alpha} ) }$ via $f_{\alpha}^*$, since 
$(f_{\alpha})_{\sharp}$ is onto.  
It follows from \S \ref{sec=init} that 
$TC_1(\V^{\alpha})=T_1(\V^{\alpha})= 
f_{\alpha}^*H^1( C_{\alpha},\C)$.
This establishes \eqref{eq=fir}.  That this 
decomposition of $\R_1(G)$ coincides with the one into irreducible 
components follows from the fact that connected algebraic 
subtori are determined by their Lie algebras; see e.g.~\cite[13.1]{Hu}.
\end{proof}

\subsection{Isotropic subspaces}
\label{subsec=isotropic}
Before proceeding, we introduce some notions which will 
be of considerable use in the sequel. 
Let $\mu \colon H^1 \wedge H^1 \to H^2$ be a $\C$-linear map, 
and $\R_k(\mu)\subset H^1$ be the corresponding resonance varieties, 
as defined in \eqref{eq=resmu}. One way to construct elements 
in these varieties is as follows.

\begin{lemma} 
\label{lem=linalg1}
Suppose $V \subset H^1$ is a linear subspace of dimension $k$. 
Set $i=\dim \im (\mu\colon V \wedge V \to H^2)$.
If $i<k-1$, then $V \subset \R_{k-i-1}(\mu)\subset \R_1(\mu)$.
\end{lemma}

\begin{proof}
Let $x \in V$, and set $x^{\perp}_V=\{y \in V \mid \mu (x \wedge y)=0\}$.  
Clearly, $\dim x^{\perp}_V \geq k - i$. On the other hand, 
$x^{\perp}_V \slash \C \cdot x \subset H^1(H^{\bullet},\mu _x)$, 
and so $x \in  \R_{k-i-1} (\mu)$. 
\end{proof}

Therefore, the subspaces $V \subset H^1$ for which 
$\dim \im (\mu\colon V \wedge V \to H^2)$ is small
are particularly interesting. This remark gives a 
preliminary motivation for the following key definition.

\begin{definition} 
\label{def=position} 
Let $\mu \colon \bigwedge ^2H^1 \to  H^2$ be a 
$\C$-linear mapping, and let $V \subset H^1$ be 
a $\C$-linear subspace. 

\begin{romenum}
\item \label{it1}
$V$ is {\em  $0$-isotropic} (or simply, {\em isotropic}) 
with respect to $\mu$ if the restriction 
$\mu ^V\colon \bigwedge ^2V \to H^2$ is trivial. 

\item \label{it2}
$V$ is  {\em $1$-isotropic}\/ with respect to $\mu$ if the restriction 
$\mu ^V\colon\bigwedge ^2V \to  H^2$ has $1$-dimensional 
image and is a non-degenerate skew-symmetric bilinear form.
\end{romenum}
\end{definition}

\begin{example} 
\label{ex=position} 
Let $C$ be a smooth curve, 
and let $\mu_C \colon \bigwedge ^2H^1(C,\C) \to H^2(C,\C)$ be 
the usual cup-product map.  There are two cases of interest to us.
\begin{romenum}
\item \label{curv1}
 If  $C$ is not compact, then $ H^2(C,\C)=0$ 
and so any subspace $V \subset H^1(C,\C)$ is isotropic.
\item \label{curv2} 
If $C$ is  compact, of genus $g\ge 1$, then 
$ H^2(C,\C)=\C$ and $ H^1(C,\C)$ is  $1$-isotropic.
\end{romenum}
\end{example}

Now let $\mu _1\colon  \bigwedge ^2H^1_1 \to  H^2_1$ and 
$\mu _2\colon  \bigwedge ^2H^1_2 \to  H^2_2$ be two $\C$-linear maps. 
\begin{definition} 
\label{def=similar} 
The maps $\mu _1$ and  $\mu _2$ are {\em equivalent} 
(notation $\mu _1 \simeq   \mu _2$) if there exist linear 
isomorphisms $\phi ^1\colon H_1^1 \to H_2^1$ and 
$\phi ^2\colon \im (\mu _1) \to \im (\mu _2)$ 
such that $\phi^2 \circ \mu _1= \mu_2 \circ \wedge ^2 \phi^1$.
\end{definition} 

The key point of this definition is that the $k$-resonant varieties 
$\R_k(\mu _1)$ and $\R_k(\mu _2)$ can be identified under $\phi ^1$ 
when $\mu _1 \simeq \mu _2$.  Moreover, subspaces that are either 
$0$-isotropic or $1$-isotropic with respect to $\mu_1$ and $\mu_2$ 
are matched under $\phi^1$.

\subsection{Admissible maps and isotropic subspaces}
\label{subsec=admiso}

We now consider in more detail which admissible maps 
$f_{\alpha}\colon M \to C_{\alpha}$ may occur in 
Theorem \ref{thm:vadm}.

\begin{prop} 
\label{prop=posobstr} 
Let $M$ be a connected quasi-compact K\"{a}hler manifold, 
and let $f\colon M \to C$ be an admissible map onto the 
smooth curve $C$.
\begin{enumerate}

\item \label{ra1} 
If $W_1(H^1(M,\C))=H^1(M,\C)$, then the curve $C$ 
is either compact, or it is obtained from a compact smooth 
curve $\oC$ by deleting a single point.

\item \label{ra2} 
If  $W_1(H^1(M,\C))=0$, then  the curve $C$ is rational.
If $\chi(C)<0$, then $C$ is obtained from $\C$ by deleting 
at least two points, and $f^* H^1( C,\C)$ is $0$-isotropic with 
respect to $\cup_M$.

\item \label{ra3} 
Assume in addition that $\pi_1(M)$ is $1$-formal.
If the curve $C$ is compact of genus at least 1, then  
$f^*\colon H^2( C,\C) \to H^2( M ,\C)$ 
is injective, and so $f^* H^1( C,\C)$ is 
$1$-isotropic with respect to $\cup_M$.

\end{enumerate}
\end{prop}

\begin{proof}
Recall that $f _{\sharp}\colon \pi_1(M) \to \pi_1(  C )$ 
is surjective; hence $f^*\colon H^1( C,\C) \to 
H^1( M,\C)$ is injective. 

Part~\eqref{ra1}.  Strictly speaking, a quasi-compact 
K\"{a}hler manifold $M$ does not have a {\em unique} mixed 
Hodge structure. Nevertheless, it inherits such a structure 
from each good compactification $\oM$, by 
Deligne's construction in the smooth quasi-projective 
case, see \cite{D}.  

By the admissibility condition on $f\colon M\to C$, there is 
a good compactification $\oM$ such that $f$ extends 
to a regular morphism ${\overline f}\colon \oM \to \oC$. 
Fixing such an extension, the condition 
$W_1(H^1(M,\C))=H^1(M,\C)$ simply means that
$j^*(H^1(\oM,\C))=H^1(M,\C)$, where $j\colon 
M \to \oM$ is the inclusion. Since regular maps 
$f$ which extend to good compactifications of source 
and target obviously preserve weight filtrations, the 
mixed Hodge structure on $H^1(C,\C)$ must be pure 
of weight $1$, see \cite{D}. If we write 
$C= \oC \setminus A$, for some finite set $A$, 
then there is an exact Gysin sequence
\[
\xymatrixcolsep{12pt}
\xymatrix{0 \ar[r]& H^1(\oC, \C) \ar[r]& 
H^1({ C}, \C) \ar[r]& H^0(A,\C)(-1) 
\ar[r]&  H^2(\oC, \C) \ar[r]& 
H^2({ C}, \C) \ar[r]& 0 },
\]
see for instance \cite[p.~246]{D1}. But $H^0(A,\C)(-1)$ is pure 
of weight $2$, and so $ H^1({ C}, \C)$ is pure of weight $1$ 
if and only if $\abs{A} \leq 1$.

\smallskip
Part~\eqref{ra2}. 
By the same argument as before, we infer in this 
case that $H^1( C,\C)$ should be pure of weight $2$.
The above Gysin sequence shows that $ H^1({ C}, \C)$ 
is pure of weight $2$ if and only if $g(\oC)=0$, 
i.e., $\oC=\PP^1$. Finally,  $\chi (C)<0$ implies 
$\abs{A} \ge 3$. 

\smallskip
Part~\eqref{ra3}. Set $G:=\pi_1(M)$, $\T_M:=\T_G$ and 
$\T:=\T_{\pi_1(C)}$. Note  that $\dim \T>0$. 
Furthermore, the character torus $\T$ is embedded in 
$\T_M$, and its Lie algebra $T_1 (\T)$ is embedded in 
$T_1(\T_M)$, via the natural maps induced by $f$.  
By Theorem \ref{thm:vadm}\eqref{va2}, 
\[
\dim H^1(M,{}_{f^*{\rho}}\C)= \dim H^1( C, {}_{\rho}\C) ,
\]
for $\rho \in \T$ near $1$ and different from $1$, since
both the surjectivity of $f_{\sharp}$ in Part \eqref{ra1}, 
and the property from Part \eqref{ra2} do not require 
the assumption $\chi (C)<0$.

Applying Theorem \ref{thm=analiso} to both $G$ (using our 
$1$-formality hypothesis), and $\pi_1(C)$ (using Example 
\ref{ex:w2formal}), we obtain from the above equality that
\[
\dim H^1(H^{\bullet}(M,\C), \mu_{f^* z})= 
\dim H^1(H^{\bullet} (C,\C), \mu_z) ,
\]
for all $z\in H^1(C, \C)$ near $0$ and different from $0$.
Moreover, for any such $z$, a standard calculation shows   
$\dim H^1(H^{\bullet} (C,\C), \mu_z)=2g-2$, 
where $g=g(C)$. 

Now suppose $f^*\colon H^2( C,\C) \to H^2( M ,\C)$ 
were not injective. Then  
$f^* H^1 ( C,\C)$ would be a $0$-isotropic 
subspace of $H^1(M, \C)$, containing $f^* (z)$.  In turn, this 
would imply  $\dim H^1(H^{\bullet}(M,\C), \mu_{f^* z})\ge 2g-1$, 
a contradiction.
\end{proof}

We close this section by pointing out the subtlety of the 
injectivity property from Proposition \ref{prop=posobstr}\eqref{ra3}.

\begin{example} 
\label{ex:notinj}
Let $L_g$ be the complex algebraic line bundle associated 
to the divisor given by a point on a projective smooth complex 
curve $C_g$ of genus $g\ge 1$. Denote by $M_g$ the total 
space of the $\C^*$-bundle associated to $L_g$. Clearly, 
$M_g$ is a smooth, quasi-projective manifold. 
(For $g=1$, this example was examined by
Morgan in \cite[p.~203]{M}.) Denote by $f_g \colon M_g\to C_g$ 
the natural projection.  This map is a locally trivial fibration, 
which is admissible in the sense of Arapura \cite{A}.
Since the Chern class $c_1(L_g)\in H^2(C_g, \Z)$ equals 
the fundamental class, it follows that
$f_g^* \colon H^2(C_g, \C)\to H^2(M_g, \C)$ is the zero map. 
Set $G_g= \pi_1(M_g)$.

A straightforward analysis of the Serre spectral sequence associated 
to $f_g$, with arbitrary untwisted field coefficients, shows that 
$(f_g)_* \colon H_1(M_g, \Z)\to H_1(C_g, \Z)$ is an isomorphism, 
which identifies the respective character tori, to be denoted in the 
sequel by $\T_g$. This also implies that $W_1(H^1 (M_g,\C))= 
H^1( M_g,\C)$, since this property holds for the compact variety $C_g$. 

We claim $f_g$ induces an isomorphism
\begin{equation} 
\label{eq:mgcg}
H^1(C_g, {}_{\rho}\C) \xrightarrow{\,\simeq\,} 
H^1(M_g, {}_{f_g^* {\rho}}\C) ,
\end{equation}
for all $\rho \in \T_g$. If $\rho=1$, this is clear. If $\rho \ne 1$, 
then $\Hom_{\Z \pi_1(C_g)}(\Z, {}_{\rho}\C)=0$, since the 
monodromy action of $\pi_1(C_g)$ on $\Z= H_1 (\C^*, \Z)$ 
is trivial. The claim follows from the $5$-term exact sequence 
for twisted cohomology associated to the group extension
$1\to \Z \to G_g \to \pi_1(C_g)\to 1$; 
see \cite[VI.8(8.2)]{HS}. 

It follows that
\begin{equation} 
\label{eq:vmg}
\V_k(G_g) = \begin{cases}
\T_g , & {\rm for}\ 0\le k\le 2g-2 ;\\
\{ 1\}  , & {\rm for}\   2g-1\le k\le 2g .
\end{cases}
\end{equation}
On the other hand, $\cup_{G_g}=0$, since $f_g^*=0$ on $H^2$. 
Therefore
\begin{equation} 
\label{eq:rmg}
\R_k(G_g) = \begin{cases}
T_1(\T_g), & {\rm for}\  0\le k\le 2g-1 ; \\
\{ 0\} ,        & {\rm for}\ k= 2g .
\end{cases}
\end{equation}

By inspecting \eqref{eq:vmg} and \eqref{eq:rmg}, 
we see that the tangent cone formula fails for $k=2g-1$. 
Consequently, the (quasi-projective) group $G_g$ 
cannot be $1$-formal. 
We thus see that the $1$-formality hypothesis from 
Proposition \ref{prop=posobstr}\eqref{ra3} is essential for 
obtaining the injectivity property of $f^*$ on $H^2$.
\end{example}

\begin{remark} 
\label{rk=posobstr} 
It is easy to show that $f^*\colon H^2( C,\C) \to H^2({M},\C)$ 
is injective when $M$ is compact. On the other hand, consider the following 
genus zero example, kindly provided to us by Morihiko Saito. Take 
$ C=\PP^1$ and  $M=\PP^1 \times \PP^1 \setminus (C_1 \cup C_2)$, 
where $C_1= \{\infty \}   \times  \PP^1 $ and $C_2$ is the diagonal in 
$\PP^1 \times \PP^1 $. The projection of $M$ on the first factor has as 
image $\C=\PP^1 \setminus \{\infty \}$ and affine lines as fibers; 
thus, $M$ is contractible. If we take $f\colon M \to \PP^1$ to be the map 
induced by the second projection, we get an admissible map such that 
$f^*\colon H^2( C,\C) \to H^2(M ,\C)$ is not injective. 
\end{remark}

\section{Position and resonance obstructions} 
\label{sec=posobs}

In this section, we give our obstructions to realizing a 
finitely presented group $G$ as the fundamental group of a 
connected, quasi-compact K\"{a}hler manifold $M$. 
We start with a definition. 

\begin{definition}
\label{def=posob} 
Let $\mu \colon H^1 \wedge H^1 \to H^2$ be a $\C$-linear map,
with $\dim H^1 >0$.  
Let $\R_k(\mu)\subset H^1$ be the corresponding resonance varieties, 
and let $\R_1(\mu)=\bigcup\nolimits_{\alpha}\R^{\alpha}$
be the decomposition of  $\R_1(\mu)$ into irreducible components. 
We say $\mu$ satisfies the {\em resonance obstructions}\/ if the 
following conditions hold. 

\begin{description}
\item[1. Linearity]  
\label{p0}  
Each component $\R^{\alpha}$ is a linear subspace of $H^1$.
\\[-8pt]
\item[2. Isotropicity] 
\label{p1}  
If $\R^{\alpha}\ne \{0\}$, then $\R^{\alpha}$ 
is a $p$-isotropic subspace of dimension 
at least $2p+2$, for some $p=p(\alpha) \in \{0,1\}$. 
\\[-8pt]
\item[3. Genericity] 
\label{p2} 
If $\alpha \ne \beta $, then $\R^{\alpha} \bigcap \R^{\beta}=\{0\}.$
\\[-8pt]
\item[4. Filtration by dimension] 
\label{p3} 
For $1\leq k\leq \dim H^1$,
\[
\R_k(\mu)=\bigcup\nolimits_\alpha \R^{\alpha},
\]
where the union is taken over all components  $\R^{\alpha}$ 
such that $\dim \R^{\alpha} >k+p(\alpha)$.  By convention, 
the union equals $\{0\}$ if the  set of such components is empty.
\end{description}
Let $\TT =\bigcup_\alpha \TT^{\alpha}$ be the irreducible
decomposition of a homogeneous variety $\TT\subset H^1$. 
We say $\TT$ satisfies the {\em position obstructions}\/ with 
respect to $\mu$ if conditions {\bf 1}--{\bf 3} above hold, 
with $\TT$ replacing $\R_1(\mu)$.
\end{definition}

Note that the above resonance conditions depend only on the 
equivalence class of $\mu$, in the sense of Definition \ref{def=similar}. 
Note also that $\R_1(\mu)=\emptyset$, if $H^1=0$, while 
$0\in \R_1(\mu)$, if $\dim H^1 >0$. If $\mu=\cup_G$, where 
$G=\pi_1(M)$, we will examine the above position conditions 
for $\TT= TC_1(\V_1(G))$. Note that, again, $\TT=\emptyset$, 
if $H^1=0$, while $0\in \TT$, if $H^1\ne 0$.

The next three lemmas will be used in establishing the position 
obstruction from Theorem \ref{thm=posobstr}\eqref{a2}.

\begin{lemma}
\label{lem=fcok} 
Let $X$ be a connected quasi-compact K\"{a}hler manifold, $C$ a smooth 
curve and $f \colon X\to C$ a non-constant holomorphic mapping. 
Assume that $f$ admits a holomorphic extension 
$\hat f\colon \hat X \to \hat C$, where $\hat X$ (resp.~$\hat C$) 
is a smooth compactification of $X$ (resp.~$C$). Then the induced 
homomorphism in homology, $f_*\colon H_1(X, \Z)\to H_1(C, \Z)$, 
has finite cokernel.
\end{lemma}

\begin{proof}
Let $Y =\Sing(\hat f)$ be the set of singular points of $\hat f$, i.e., 
the set of all points $x \in \hat X$ such that $d_x\hat f=0$. Then 
$Y$ is a closed analytic subset of $\hat X$. Using Remmert's Theorem, 
we find that $Z=\hat f (Y)$ is a closed analytic subset of $\hat C$,  
and $\hat f(\hat X)=\hat C$.  By Sard's Theorem, $Z \ne \hat C$, 
hence $Z$ is a finite set. 

Let $B=(\hat C \setminus C) \cup Z$; set $C'=\hat C \setminus B$, 
and $\hat X '=\hat X \setminus \hat f^{-1}(B)   =\hat f^{-1}(C')$. Then 
the restriction $\hat f '\colon \hat X '\to C'$ is a locally trivial fibration;  
its fiber is a compact manifold, and thus has only finitely many 
connected components. Using the tail end of the homotopy 
exact sequence of this fibration, we deduce 
that the induced homomorphism, 
$\hat f'_{\sharp} \colon \pi_1(\hat X ')\to \pi_1(C')$, 
has image of finite index. 

Now note that $i \circ \hat f ' \circ k=f\circ j$, where 
$i\colon C' \to C$, $j\colon X \setminus \hat f^{-1}(B) \to X$, 
and $k\colon X \setminus \hat f^{-1}(B) \to \hat X '$ are the 
inclusion maps. From the above, it follows that
$f_{*}\colon H_1(X,\Z)\to H_1(C,\Z)$ has image 
of finite index. 
\end{proof}

Let $\V^{\alpha}$ and $\V^{\beta}$ be two distinct, positive-dimensional, 
irreducible components of $\V_1(\pi_1(M))$ containing $1$. Realize 
them by pull-back, via admissible maps, $f_{\alpha}\colon M\to C_{\alpha}$ and 
$f_{\beta}\colon M\to C_{\beta}$, as in Theorem \ref{thm:vadm}\eqref{va1}.
We know that generically (that is, for $t\in C_{\alpha}\setminus B_{\alpha}$, 
where $B_{\alpha}$ is finite) the fiber $f_{\alpha}^{-1}(t)$  is smooth and 
irreducible.

\begin{lemma}
\label{lem=noncst}
In the above setting, there exists $t\in C_{\alpha}\setminus B_{\alpha}$
such that the restriction of $f_{\beta}$ to $f_{\alpha}^{-1}(t)$ 
is non-constant.
\end{lemma}

\begin{proof}
Assume $f_{\beta}$ has constant value, $h(t)$, on the fiber 
$f_{\alpha}^{-1}(t)$, for $t\in C_{\alpha}\setminus B_{\alpha}$. 
We first claim that this implies the existence of a continuous 
extension, $h\colon C_{\alpha}\to C_{\beta}$, with the property 
that $h\circ f_{\alpha}= f_{\beta}$. 

Indeed, let us pick an arbitrary special value, $t_0\in B_{\alpha}$, 
together with a sequence of generic values, 
$t_n\in C_{\alpha}\setminus B_{\alpha}$, 
converging to $t_0$. For any $x\in f_{\alpha}^{-1}(t_0)$, 
note that the order at $x$ of the holomorphic function $f_{\alpha}$ 
is finite. Hence, we may find a sequence, $x_n \to x$, such that 
$f_{\alpha}(x_n)=t_n$. By our assumption, $f_{\beta}(x)= \lim h(t_n)$, 
independently of $x$, which proves the claim.

At the level of character tori, the fact that 
$h\circ f_{\alpha}= f_{\beta}$ implies  
$\V^{\beta}= f_{\beta}^* \T_{\pi_1(C_{\beta})}\subset 
f_{\alpha}^*\T_{\pi_1(C_{\alpha})}= \V^{\alpha}$, a contradiction.
\end{proof}

\begin{lemma}
\label{lem=finite}
Let $\V^{\alpha}$ and $\V^{\beta}$ be two distinct 
irreducible components of $\V_1(\pi_1(M))$ containing $1$.  
Then $\V^{\alpha}\cap \V^{\beta}$ is finite.
\end{lemma}

\begin{proof}
We may suppose that both components are positive-dimensional. 
Lemma \ref{lem=noncst} guarantees the existence of a 
generic fiber of $f_{\alpha}$, say $F_{\alpha}$, with the property 
that the restriction of $f_{\beta}$ to $F_{\alpha}$, call it 
$g\colon F_{\alpha} \to C_{\beta}$, is non-constant. By 
Lemma \ref{lem=fcok}, there exists a positive integer $m$
with the property that
\begin{equation}
\label {eq=cokf}
m\cdot H_1(C_{\beta}, \Z)\subset \im (g_*)\, .
\end{equation}
We will finish the proof by showing that
$\rho^m =1$, for any $\rho\in \V^{\alpha}\cap \V^{\beta}$.

To this end, write $\rho= \rho_{\beta}\circ (f_{\beta})_*$, 
with $\rho_{\beta}\in \T_{\pi_1(C_{\beta})}$. For an 
arbitrary element $a\in H_1(M, \Z)$, we have 
$\rho^m (a)= \rho_{\beta}(m\cdot (f_{\beta})_*a)$. 
From \eqref{eq=cokf}, it follows that 
$m\cdot (f_{\beta})_*a=(f_{\beta})_*( j_\alpha)_*b$, 
for some $b\in H_1(F_{\alpha}, \Z)$, where 
$j_\alpha\colon F_{\alpha}\hookrightarrow M$ is 
the inclusion. 
On the other hand, we may also write 
$\rho= \rho_{\alpha}\circ (f_{\alpha})_*$, with 
$\rho_{\alpha}\in \T_{\pi_1(C_{\alpha})}$. Hence,  
$\rho^m (a)= \rho_{\alpha}((f_{\alpha})_* (j_{\alpha})_*b)= 
\rho_{\alpha}(0)=1$, as claimed.
\end{proof}

The next result establishes Parts \eqref{a1}--\eqref{a3} of 
Theorem \ref{thm=posobstr} from the Introduction. 

\begin{theorem} 
\label{thm=posob} 
Let $M$ be a connected, quasi-compact K\"{a}hler manifold. 
Set $G=\pi_1(M)$ and assume $b_1(G)>0$.  
Let $\{ \V^{\alpha} \}$ be the irreducible 
components of $\V_1(G)$ containing $1$.
Then the following hold.

\begin{enumerate}

\item \label{a61}
The tangent spaces $\TT^{\alpha}:= T_1(\V^{\alpha})$ are the
irreducible components of $\TT:= TC_1(\V_1(G))$.

\item \label{a62}
If $G$ is $1$-formal, $\{ \TT^{\alpha}\}$ is the family of irreducible
components of $\R_1(G)$.

\item \label{a63}
The variety $\TT$ satisfies the position obstructions 
{\bf 1}--{\bf 3} from Definition \ref{def=posob}, with respect to 
$\cup_G\colon H^1(G,\C) \wedge H^1(G,\C)\to H^2(G,\C)$.

\item \label{a64}
If $G$ is $1$-formal, $\cup_G$ verifies the resonance conditions 
{\bf 1}--{\bf 4} from \ref{def=posob}.
\end{enumerate}
\end{theorem}

\begin{proof}
Part \eqref{a61}. As noticed before, $\TT= \bigcup_{\alpha} \TT^{\alpha}$.
By \cite[13.1]{Hu}, $\TT= \bigcup_{\alpha} \TT^{\alpha}$ is the  
decomposition of $\TT$ into irreducible components.

Part \eqref{a62}. If $\R_1(G)=0$, then $\TT=0$, and there is 
nothing to prove. If $\R_1(G)\ne 0$, the statement follows 
from Corollary \ref{cor:radm}.

Part \eqref{a63}. Property {\bf 1} is clear. As for property {\bf 2}, we have 
$\V^{\alpha}= f_{\alpha}^*\T_{\pi_1(C_{\alpha})}$, where $f_{\alpha}$ is
admissible and $\chi(C_{\alpha})<0$; see Theorem \ref{thm:vadm}\eqref{va1}.
Therefore, $\TT^{\alpha}= f_{\alpha}^* H^1(C_{\alpha}, \C)$. If the curve
$C_{\alpha}$ is non-compact, the subspace $\TT^{\alpha}$ is clearly isotropic,
and $\dim \TT^{\alpha}=b_1(C_\alpha) \ge 2$. If $C_{\alpha}$ is compact and
$f_{\alpha}^*$ is zero on $H^2(C_{\alpha}, \C)$, we obtain 
the same conclusion as before. Finally, if $C_{\alpha}$ is compact and
$f_{\alpha}^*$ is non-zero on $H^2(C_{\alpha}, \C)$, then plainly
$\TT^{\alpha}$ is $1$-isotropic and $\dim \TT^{\alpha}=b_1(C_\alpha) \ge 4$.
The isotropicity property is thus established.

By \cite[Theorem 12.5]{Hu},  property {\bf 3} 
is equivalent to the fact that $\V^{\alpha}\cap \V^{\beta}$ 
is finite. Thus, the genericity property follows from 
Lemma \ref{lem=finite}.

Part \eqref{a64}. By Parts \eqref{a61}--\eqref{a63}, the resonance 
conditions {\bf 1}--{\bf 3} are verified  by $\cup_G$. Property {\bf 4} 
may be checked as follows.  We may assume $\R_1(G)\ne \{0\}$, 
since otherwise the property holds trivially.  We may also assume 
$k<b_1(G)$, since otherwise $\R_k(G)= \{0\} $, and there is nothing 
to prove.

To prove the desired equality, we have to check first that any 
non-zero element $u \in \R_k(G)$ belongs to some 
$ \R^{\alpha}$ with $\dim \R^{\alpha} >k$. 
Definition  \eqref{eq=resmu} guarantees the existence of 
elements $v_1,\dots ,v_k \in H^1(G,\C)$ with $v_i \cup  u=0$, 
and such that $u, v_1,\dots ,v_k$ are linearly independent. 
Since the subspaces  $\langle u, v_i \rangle $ spanned by the 
pairs $\{u,v_i\}$ are clearly contained in $\R_1(G)$, it follows that 
$\langle u,v_i \rangle\subset  \R^{\alpha _i}$.  
Necessarily $ \alpha _1=\cdots =\alpha _k :=\alpha$, since 
otherwise property {\bf 3} would be violated. This proves 
$u\in \R^{\alpha}$, with $\dim  \R^{\alpha} >k$. 

Now, if $p({\alpha})=1$, then $\dim \R^{\alpha}>k+1$.  
For otherwise, we would have $\R^{\alpha}= u^{\perp}_{ \R^{\alpha}}$, 
which would violate the non-degeneracy property from 
Definition \ref{def=position}\eqref{it2}.
Finally, that $\dim  \R^{\alpha} >k+ p({\alpha})$ implies 
$\R^{\alpha} \subset  \R_k(G)$ follows at once  
from Lemma \ref{lem=linalg1}. 
\end{proof}

We now relate the dimensions of 
the cohomology groups $H^1(H^{\bullet} (M,\C), \mu_z)$ 
and $H^1(M, {}_{\rho}\C)$ corresponding to 
$z\in \Hom(G,\C)\setminus \{0\}$ and 
$\rho=\exp(z)\in \Hom(G,\C^*)\setminus \{1\}$, 
to the dimension and isotropicity of the resonance 
component $\R^{\alpha}$ to which $z$ belongs.

\begin{prop}
\label{prop=p4} 
Let $M$ be a connected  quasi-compact K\"{a}hler manifold, 
with fundamental group $G$, and first resonance variety 
$\R_1(G)=\bigcup_{\alpha}  \R^{\alpha}$.  If $G$ is $1$-formal, 
then the following hold. 

\begin{enumerate}
\item \label{f1}
If $z\in \R^{\alpha}$ and $z\ne 0$, then 
$\dim H^1(H^{\bullet} (M,\C), \mu_z)=\dim \R^{\alpha}-p(\alpha)-1$.

\item \label{f2}
If $\rho \in \exp (\R^{\alpha})$ and $\rho\ne 1$, then 
$\dim H^1(M, {}_{\rho}\C)\ge \dim \R^{\alpha}-p(\alpha)-1$, with 
equality for all except finitely many local systems $\rho$.
\end{enumerate}
\end{prop}

\begin{proof} 
Part~\eqref{f1}. 
Recall that $\R^{\alpha}=f_{\alpha}^* H^1(C_{\alpha}, \C)$. 
Exactly as in the proof of Proposition \ref{prop=posobstr}\eqref{ra3}, 
we infer that 
\begin{equation} 
\label{eq:dimgen}
\dim H^1(M, {}_{\rho}\C)= 
\dim H^1(H^{\bullet}(M,\C), \mu_z)= 
\dim H^1(H^{\bullet} (C_{\alpha},\C), \mu_{\zeta}) ,
\end{equation}
where $z=f_{\alpha}^* \zeta$ and $\rho =\exp (z)$,
for all $z\in \R^{\alpha}$ near $0$ and different from $0$. Clearly
\begin{equation} 
\label{eq:dimcurve}
\dim H^1(H^{\bullet} (C_{\alpha},\C), \mu_{\zeta})= 
\dim \R^{\alpha}- p({\alpha})-1 ,
\end{equation}
if $\zeta\ne 0$. Since plainly 
$\dim H^1(H^{\bullet}(M,\C), \mu_z)= 
\dim H^1(H^{\bullet}(M,\C), \mu_{\lambda z})$, 
for all $\lambda \in \C^*$,
equations \eqref{eq:dimgen} and \eqref{eq:dimcurve} 
finish the proof of \eqref{f1}.

Part~\eqref{f2}. 
Starting from the standard presentation of the group 
$\pi_1(C_{\alpha})$, a Fox calculus computation shows that 
$\dim H^1(C_{\alpha}, {}_{\rho'} \C)=\dim \R^{\alpha}- p({\alpha})-1$, 
provided $\rho' \ne 1$. By Theorem \ref{thm:vadm}\eqref{va2}, the 
equality $\dim H^1(M, {}_{\rho}\C)= \dim \R^{\alpha}- p({\alpha})-1$ 
holds for all but finitely many local systems $\rho\in \exp(\R^{\alpha})$. 
By semi-continuity (see Remark \ref{rk=semi}), the inequality
$\dim H^1(M, {}_{\rho}\C)\ge \dim \R^{\alpha}- p({\alpha})-1$ 
holds for all $\rho \in \exp(\R^{\alpha})$. 
\end{proof}

Using Proposition \ref{prop=posobstr}, we obtain the 
following corollary to Theorem \ref{thm=posob}. 

\begin{corollary} 
\label{cor:firstex}
Let $M$ be a connected  quasi-compact K\"{a}hler manifold, 
with fundamental group $G$, and first resonance variety 
$\R_1(G)=\bigcup_{\alpha}  \R^{\alpha}$. Assume $b_1(G)>0$ and
$\R_1(G)\ne \{0\}$.

\begin{enumerate}
\item \label{appl1} 
If $M$ is compact then $G$ is $1$-formal, and $\cup_G$ satisfies 
the resonance obstructions.  Moreover, each component 
$\R^{\alpha}$ is $1$-isotropic, with $\dim \R^{\alpha}= 2g_{\alpha}\ge 4$.

\item \label{appl2} 
If $W_1(H^1 (M,\C))=0$ then $G$ is $1$-formal, and $\cup_G$ satisfies 
the resonance obstructions.  Moreover, each component 
$\R^{\alpha}$  is  $0$-isotropic, with $\dim \R^{\alpha}\ge 2$.

\item \label{appl3} 
If $W_1(H^1 (M,\C))=H^1 (M,\C)$ and $G$ is $1$-formal, then
$\dim \R^{\alpha}= 2g_{\alpha}\ge 2$, for all $\alpha$.

\end{enumerate}
\end{corollary}

Next, we sharpen a result of Arapura \cite[Corollary V.1.9]{A}, 
under a $1$-formality assumption, thereby extending known 
characterizations of  `fibered' K\"{a}hler groups 
(see \cite[\S2.3]{ABC} for a survey).

\begin{corollary}
\label{cor:fibered}
Let $M$ be a connected  quasi-compact K\"{a}hler manifold. 
Suppose the group $G=\pi_1(M)$ is $1$-formal. 
The following are then equivalent.

\begin{romenum}
\item \label{gf1}
There is an admissible map,
$f \colon M\to C$, onto a smooth complex curve $C$ with 
$\chi (C)<0$.
\item \label{gf2}
There is an epimorphism, $\varphi\colon G\surj \pi_1(C)$, onto 
the fundamental group of a 
smooth complex curve $C$ with $\chi (C)<0$.
\item \label{gf3}
There is an epimorphism, $\varphi\colon G\surj \bF_r$, onto
a free group of rank $r\ge 2$.
\item \label{gf4}
The resonance variety $\R_1(G)$ strictly contains $\{ 0\}$.
\end{romenum}
\end{corollary}

\begin{proof}
The equivalence of the first three properties is proved in 
\cite[Corollary V.1.9]{A}, without our additional $1$-formality 
hypothesis. The implication \eqref{gf3} $\Rightarrow$
\eqref{gf4} also holds in general; this may be seen by 
using the $\cup_G$-isotropic subspace 
$\varphi^* H^1(\bF_r, \C)\subset H^1(G, \C)$.
Finally, assuming $G$ is $1$-formal, the implication 
\eqref{gf4} $\Rightarrow$ \eqref{gf1} follows from 
Corollary \ref{cor:radm}. 
\end{proof}

The equivalence \eqref{gf3} $\Leftrightarrow$ \eqref{gf4} 
above finishes the proof of Theorem \ref{thm=posobstr} 
from the Introduction. 

We close this section with a pair of examples showing that 
both the quasi-K\"{a}hler  and the $1$-formality 
assumptions are needed in order for this equivalence 
to hold.  

\begin{example}
\label{ex:heisenberg}
Consider the smooth, quasi-projective variety $M_1$ 
from Example \ref{ex:notinj} (the complex version of the 
Heisenberg manifold). As mentioned previously, the 
group $G_1=\pi_1(M_1)$ is not $1$-formal. On the 
other hand, it is well-known that $G_1$ is a nilpotent 
group.   Therefore, property \eqref{gf3} fails for $G_1$. 
Nevertheless, $\R_1(G_1)=\C^2$, by  virtue of 
\eqref{eq:rmg}, and so property \eqref{gf4} holds
for $G_1$.
\end{example}

\begin{example}
\label{ex:nonorientable}
Let $N_h$ be the non-orientable  surface of genus 
$h\ge 1$, that is, the connected sum of $h$ real projective planes. 
It is readily seen that $N_h$ has the rational homotopy type of 
a wedge of $h-1$ circles.  Hence $N_h$ is a formal space, 
and so $\pi_1(N_h)$ is a $1$-formal group.   
Moreover, $\R_1(\pi_1(N_h))=\R_1(\bF_{h-1})$, and so 
$\R_1(\pi_1(N_h))=\C^{h-1}$, provided $h\ge 3$.  
Thus, property \eqref{gf4} holds for all groups 
$\pi_1(N_h)$ with $h\ge 3$. 

Now suppose there is an epimorphism 
$\varphi\colon\pi_1(N_h) \surj \bF_r$ with $r\ge 2$, as in \eqref{gf3}.  
Then the subspace 
$\varphi^* H^1(\bF_r, \Z_2)\subset H^1(\pi_1(N_h), \Z_2)$ 
has dimension at least $2$, and is isotropic with respect to 
$\cup_{\pi_1(N_h)}$.  Hence, $h\ge 4$, by Poincar\'e duality 
with $\Z_2$ coefficients.

Focussing on the case $h=3$, we see that the group $\pi_1(N_3)$ 
is $1$-formal, yet the implication \eqref{gf4} $\Rightarrow$ \eqref{gf3} 
from Corollary \ref{cor:fibered} fails for  this group.
It follows that $\pi_1(N_3)$ cannot be realized as the fundamental 
group of a  quasi-compact K\"{a}hler manifold. Note 
that this assertion is {\em not}\/ a consequence of Theorem
\ref{thm=posob}\eqref{a64}; indeed, 
$\cup_{\pi_1(N_h)} \simeq \cup_{\bF_{h-1}}$ (over $\C$), 
while $\bF_{h-1}=\pi_1(\PP^1 \setminus \{\text{$h$ points}\})$, 
for all $h\ge 1$.
\end{example}

\section{Wedges and products} 
\label{sec=wp}

In this section, we analyze products and coproducts of groups, 
together with their counterparts at the level of first 
resonance varieties.  Using our obstructions, 
we obtain conditions for realizability of free products 
of  groups by  quasi-compact K\"{a}hler manifolds. 

\subsection{Products, coproducts, and $1$-formality}
\label{subsec=cp1}

Let $\bF(X)$ be the free group on a finite set $X$, and let 
$ {\bL}^*(X)$ be the free Lie algebra on $X$, over a field $\K$ 
of characteristic $0$.  Denote by $\widehat {\bL}(X)$ the 
Malcev Lie algebra obtained from $ {\bL}^*(X)$ by completion 
with respect to the degree filtration. Define the group 
homomorphism $\kappa _X\colon \bF(X) \to 
\eexp(\widehat {\bL}(X))$ by $\kappa _X(x)=x$ for $x \in X$. 
Standard commutator calculus \cite{L} shows that
\begin{equation} 
\label{eq=malfree}
\gr^*(\kappa _X)\colon \gr^*(\bF(X)) \otimes \K  \xrightarrow{\,\simeq\,} 
\gr^*_F(\widehat {\bL}(X))
\end{equation}
is an isomorphism. 
It follows from \cite[Appendix A]{Q} that $\kappa _X$ is a Malcev completion. 

Now let $G$ be a finitely presented group, with presentation
$G= \langle x_1,\dots ,x_s  \mid w_1,\dots ,w_r \rangle$, or, 
for short, $G=\bF(X)/\langle { \bf w} \rangle$.  Denote by 
$\llangle{ \bf w}\rrangle$ 
the closed Lie ideal of $\widehat {\bL}(X)$ generated by 
$\kappa _X(w_1), \dots ,$ $ \kappa _X(w_r)$, and consider 
the group morphism induced by $\kappa _X$,
\begin{equation} 
\label{eq=malgen}
\kappa _G\colon G \to 
\eexp \big(\widehat {\bL}(X)/\llangle{ \bf w}\rrangle\big) .
\end{equation}
It follows from \cite{P1} that $\kappa_G$ is a Malcev completion 
for $G$. (For the purposes of that paper, it was assumed that 
$G_{\ab}$ had no torsion, see  \cite[Example 2.1]{P1}. 
Actually, the proof of the Malcev completion property 
applies verbatim to the general case, 
see \cite[Theorem 2.2]{P1}.)

\begin{prop} 
\label{prop=products}
If $G_1$ and $G_2$ are $1$-formal groups, then their 
coproduct $G_1 *G_2$ and their product $G_1 \times G_2$ 
are again $1$-formal groups.
\end{prop}

\begin{proof}
First consider two arbitrary finitely presented groups, 
with presentations  $G_1= \bF(X)/\langle{ \bf u}\rangle$ 
and $G_2= \bF(Y)/\langle{ \bf v}\rangle $.  Then 
$G_1 *G_2=\bF(X \cup Y)/\langle {\bf u}, {\bf v} \rangle $. 
It follows from \eqref{eq=malgen} that  
$E_{G_1 *G_2}=E_{G_1 } \coprod   E_{G_2}$, 
the coproduct Malcev Lie algebra. 

On the other hand, $G_1 \times G_2=\bF(X \cup Y)/
\langle{ \bf u}, { \bf v}, (x,y); x \in X, y \in Y \rangle$, 
and so, by the same reasoning, 
$E_{G_1 \times G_2}= 
\widehat {\bL}(X \cup Y)/\llangle\kappa_X({ \bf u}), 
\kappa _Y({ \bf v}), (x,y); x \in X, y \in Y  \rrangle$.  
Using the Campbell-Hausdorff formula, we may 
replace each CH-group commutator $(x,y)$ with the 
corresponding Lie bracket, $[x,y]$; see 
\cite[Lemma 2.5]{P2} for details. We conclude that
$E_{G_1 \times G_2}= E_{G_1 } \prod   E_{G_2}$, 
the product Malcev Lie algebra.

Now assume  $G_1$ and $G_2$ are $1$-formal.  
In view of  Lemma \ref{lem=quadr}, we may write 
$E_{G_1 }= \widehat {\bL}(X')/\llangle{ \bf u'}\rrangle $
and $E_{G_2 }= \widehat {\bL}(Y')/\llangle{ \bf v'}\rrangle$, 
where the defining relations ${ \bf u'}$ and ${ \bf v'}$ are 
quadratic. Hence
\begin{align*} 
\label{eq=coprod}
E_{G_1* G_2}& =\widehat {\bL}(X'\cup Y')/
\llangle{ \bf u'},{ \bf v'} \rrangle ,
\\
E_{G_1 \times G_2}&= \widehat {\bL}(X'\cup Y')/
\llangle{ \bf u'},{ \bf v'},[x',y']; x' \in X', y' \in Y' \rrangle  .
\end{align*}
Since the relations in these presentations are clearly quadratic, 
the $1$-formality of both $G_1* G_2$ and $G_1 \times G_2$ 
follows from Lemma \ref{lem=quadr}.  
\end{proof}

\subsection{Products, coproducts, and resonance}
\label{subsec=cpr}

Let $U^{i}$, $V^{i}$ ($i=1,2$) be complex vector spaces.  
Given two $\C$-linear maps, 
$\mu _U\colon U^1 \wedge  U^1 \to U^2$ 
and  $\mu _V\colon V^1 \wedge  V^1 \to V^2$, 
set $W^i=U^i \oplus V^i$, and define 
$\mu _U * \mu _V\colon W^1 \wedge  W^1 \to W^2$ 
as follows:  
\[
{\mu _U * \mu _V}\!\left|_{ U^1 \wedge  U^1} \right. = \mu_U, 
\quad
{\mu _U * \mu _V}\!\left|_{ V^1 \wedge  V^1} \right. = \mu_V, 
\quad
{\mu _U * \mu _V}\!\left|_{ U^1 \wedge  V^1} \right. = 0.
\]
When $\mu _U=\cup _{G_1}$ and $\mu _V=\cup _{G_2}$, 
then clearly $\mu _U * \mu _V=\cup _{G_1*G_2}$, since 
$K(G_1*G_2,1)=K(G_1,1) \vee K(G_2,1)$.

\begin{lemma} 
\label{lem=reswedge}
Suppose $U^1 \ne 0$,  $V^1 \ne 0$, and  $\mu_U * \mu_V$ 
satisfies the second resonance obstruction from 
Definition \ref{def=posob}.  Then $\mu_U= \mu_V=0$.
\end{lemma}

\begin{proof}
Set $\mu:=\mu _U * \mu _V$.  
We know $\R_1(\mu)=W^1$, by \cite[Lemma 5.2]{PS1}. 
If $\mu \ne 0$, then
$\mu$ is $1$-isotropic, with  $1$-dimensional 
image. It follows that either  $\mu _U= 0$ or $\mu _V=0$.
In either case, $\mu$ fails to be non-degenerate, a contradiction. 
Thus, $\mu=0$, and so $\mu _U= \mu _V=0$.
\end{proof}

Next, given $\mu_U$ and  $\mu _V$ as above, set 
$Z^1=U^1 \oplus V^1$ and 
$Z^2=U^2 \oplus V^2 \oplus (U^1 \otimes V^1)$, 
and define  
$\mu _U \times \mu _V\colon Z^1 \wedge  Z^1 \to Z^2$ 
as follows. As before, the restrictions of $\mu _U \times \mu _V$ to 
$U^1 \wedge  U^1$ and $V^1 \wedge  V^1$ are given by 
$\mu _U$ and  $\mu _V$,  respectively. On the other hand,  
$\mu _U \times \mu _V(u \wedge v)=u \otimes v$, for 
$u \in U^1$ and $v \in  V^1$.  
Finally, if $\mu _U=\cup _{G_1}$ and $\mu _V=\cup _{G_2}$, 
then $\mu _U \times \mu _V=\cup _{G_1 \times G_2}$, since 
$K(G_1\times G_2,1)=K(G_1,1) \times K(G_2,1)$.

\begin{lemma} 
\label{lem=resprod}
With notation as above, 
$\R_1(\mu _U \times \mu _V)=\R_1(\mu _U)\times \{0\} 
\cup\{0\}\times  \R_1(\mu _V)$.
\end{lemma}

\begin{proof}
Set $\mu=\mu_U \times \mu_V$.  The inclusion $\R_1(\mu)
\supset \R_1(\mu _U)\times \{0\} \cup\{0\}\times  \R_1(\mu _V)$ 
is obvious.  To prove the other inclusion, assume 
$\R_1(\mu) \ne 0$ (otherwise, there is nothing to prove), 
and pick $0\ne a+b \in \R_1(\mu)$, with $a \in U^1$ and 
$b \in V^1$. By definition of $\R_1(\mu)$, there is 
$x+y \in U^1 \oplus V^1$ such that 
$(a+b) \wedge (x+y) \ne 0$ and  
\begin{equation} 
\label{eq=muvanish}
\mu ((a+b) \wedge (x+y))= \mu_U(a\wedge x)+
\mu_V(b  \wedge y)+ a \otimes y -x  \otimes b=0 .
\end{equation}
In particular, $a \otimes y =x  \otimes b$. There are several 
cases to consider.

If $a \ne 0$ and $b\ne 0$, we must have $x=\lambda a$ 
and $y=\lambda b$, for some $\lambda \in \C$, and so 
$(a+b) \wedge (x+y) =(a+b)\wedge \lambda (a+b)=0$, 
a contradiction. 

If $b =0$, then $a \ne 0$ and \eqref{eq=muvanish}
forces $y=0$ and $\mu_U(a  \wedge x)=0$. 
Since $(a+b) \wedge (x+y)=a \wedge x \ne 0$, it follows
that $a \in \R_1(\mu _U)$, as needed. The other case, $a=0$, 
leads by the same reasoning to $b \in  \R_1(\mu _V)$.
\end{proof}

If $G_1$ and  $G_2$ are finitely presented groups, 
Lemma \ref{lem=resprod} implies that 
$\R_1(G_1 \times G_2)= \R_1(G_1) \times \{0\} \cup \{0\}\times \R_1(G_2)$.
An analogous formula holds for the characteristic varieties: 
$\V_1(G_1 \times G_2)=\V_1(G_1) \times \{1\}
\cup \{1\}\times \V_1(G_2)$, see \cite[Theorem 3.2]{CS2}. 

\subsection{Quasi-projectivity of coproducts} 
\label{subsec=vs}

Here is an application of Theorem \ref{thm=posobstr}. 
It is inspired by a result of M.~Gromov, who proved in 
\cite{G} that no non-trivial free product of groups can 
be realized as the fundamental group of a compact 
K\"{a}hler manifold.  We need two lemmas.

\begin{lemma} 
\label{lem=cup0}
Let $G$ be a finitely presented, commutator relators group
(that is, $G=\bF(X)/\langle {\bf w} \rangle $, with $X$ and 
${\bf w}$ finite, and ${\bf w}\subset \Gamma_2 \bF(X)$). 
Suppose $G$ is $1$-formal, and $\cup _{G}=0$.  
Then $G$ is a free group. 
\end{lemma} 

\begin{proof}
Pick a presentation $G=\bF(X)/\langle {\bf w} \rangle $, with 
all relators $w_i$ words in the commutators $(g,h)$, where 
$g,h \in \bF(X)$. We have 
$E_G=\widehat{\h} (G)$, by the $1$-formality of $G$, 
and $\h(G)=\bL (X)$, by the vanishing of $\cup _{G}$.  
Hence, $E_G= \widehat {\bL}(X)$. 
On the other hand, \eqref{eq=malgen} implies 
$E_G= \widehat {\bL}(X)/\llangle{\bf w} \rrangle$. 
We thus obtain a filtered Lie algebra isomorphism,  
$\widehat {\bL}(X) \xrightarrow{\,\simeq\,} 
\widehat {\bL}(X)/\llangle{\bf w} \rrangle$.

Taking quotients relative to the respective Malcev filtrations 
and comparing vector space dimensions, we see that 
$\kappa_X (w_i) \in \bigcap_{k \geq 1}F_k\widehat {\bL}(X)=0$, 
for all $i$. A well-known result of Magnus (see \cite{MKS}) says 
that $\gr^*(\bF(X))$ is a torsion-free graded abelian group. 
We infer from \eqref{eq=malfree} that 
$w_i \in  \bigcap_{k \geq 1} \Gamma_k\bF(X)$, for all $i$. 
Another well-known result of Magnus  (see \cite{MKS})  insures 
that $\bF(X)$ is residually nilpotent, i.e., 
$\bigcap_{k \geq 1} \Gamma_k\bF(X)=1$. 
Hence, $w_i=1$,  for all $i$, and so $G= \bF(X)$.
\end{proof}

\begin{lemma}
\label{lem=vfox}
Let $G_1$ and  $G_2$ be finitely presented groups with 
non-zero first Betti number. Then 
$\V_1(G_1 * G_2)=\T_{G_1 * G_2}$.
\end{lemma}

\begin{proof}
Let $G= \langle x_1, \dots, x_s \mid  w_1, \dots, w_r \rangle$ be 
an arbitrary finitely presented group, and let $\rho\in \T_G$ be an
arbitrary character. Recall from Remarks \ref{rk=semi} and 
\ref{rem=compu} that $\rho\in \V_1(G)$ if and only if 
$b_1(G, \rho)>0$, where
$b_1(G, \rho) := \dim \ker d_1(\rho)- \rank d_2(\rho)$. 
Moreover, the linear map $d_1(\rho)\colon \C^s \to \C$ 
sends the basis element corresponding to the generator 
$x_i$ to $\rho (x_i)-1$, while the linear map 
$d_2(\rho)\colon \C^r \to \C^s$ is given by the 
evaluation at $\rho$ of the matrix of free derivatives 
of the relators, 
$\big( \frac{\partial w_j}{\partial x_i} (\rho) \big)$; 
see Fox \cite{Fox}.

For $G=G_1 * G_2$, write $\rho= (\rho_1, \rho_2)$, 
with $\rho_i \in \T_{G_i}$. We then have 
$d_j(\rho)= d_j(\rho_1) +d_j(\rho_2)$,  for $j=1, 2$. Hence, 
$b_1(G, \rho)= b_1(G_1, \rho_1) +  b_1(G_2, \rho_2) +1$, if both
$\rho_1$ and $\rho_2$ are different from $1$, and otherwise
$b_1(G, \rho)= b_1(G_1, \rho_1) +  b_1(G_2, \rho_2)$. Since
$b_1(G_i, 1)= b_1(G_i)>0$, the claim follows. 
\end{proof}

\begin{theorem} 
\label{thm=nonalgw}
Let $G_1$ and  $G_2$ be finitely presented groups 
with non-zero first Betti number.

\begin{enumerate}
\item \label{w1}
If the coproduct $G_1 * G_2$ is quasi-K\"{a}hler, then 
$\cup_{G_1}=\cup_{G_2}=0$.

\item \label{w2}
Assume moreover that $G_1$ and  $G_2$ are $1$-formal, 
presented by commutator relators only. Then 
$G_1 * G_2$ is a quasi-K\"{a}hler group if and only
if both $G_1$ and  $G_2$ are free.
\end{enumerate}
\end{theorem}

\begin{proof}
Part \eqref{w1}. Set $G= G_1 * G_2$. From Lemma \ref{lem=vfox}, 
we know that there is just one irreducible component of $\V_1(G)$ 
containing $1$, namely $\V =\T_G^{0}$, the component of the identity
in the character torus. Hence, $T_1(\V)= H^1(G, \C)$. Libgober's 
result from \cite{Li} implies then that $\R_1(G)=H^1(G, \C)$. 
If $G$ is quasi-K\"{a}hler, Theorem \ref{thm=posobstr}\eqref{a1} 
may be invoked to infer that $\cup_G$ satisfies the isotropicity 
resonance obstruction from Definition \ref{def=posob}. The 
conclusion follows from Lemma \ref{lem=reswedge}.

Part \eqref{w2}. 
If $G_1$ and $G_2$ are free, then $G_1 * G_2$ is also free 
(of finite rank), thus quasi-projective. For the converse, 
use Part \eqref{w1} to deduce that $\cup _{G_1}=\cup _{G_2}=0$, 
and then apply Lemma \ref{lem=cup0}. 
\end{proof}

Let $\CC$ be the class of fundamental groups of complex 
projective curves of non-zero genus. Each $G \in \CC$ is a 
$1$-formal group, admitting a presentation with a single 
commutator relator, and is not free (for instance, since 
$\cup _G \ne 0$).  Proposition \ref{prop=products} and 
Theorem \ref{thm=nonalgw} yield the following corollary.  

\begin{corollary} 
\label{cor=cop1}
If $G_1,G_2 \in \CC$, then $G_1 *G_2$ is 
a $1$-formal group, yet $G_1 *G_2$ is not realizable as 
the fundamental group of a smooth, quasi-projective variety $M$. 
\end{corollary}

This shows that $1$-formality and quasi-projectivity may 
exhibit contrasting behavior with respect to the coproduct 
operation for groups.

\section{Arrangements of real planes}
\label{sec=realarr}

Let $\A=\{H_1,\dots,H_n\}$ be an arrangement of planes 
in $\RR^4$, meeting transversely at the origin.  By 
intersecting $\A$ with a $3$-sphere about $0$, we 
obtain a link $L$ of $n$ great circles in $S^3$. 
It is readily seen that the complement $M$ of the 
arrangement deform-retracts onto the complement 
of  the link.  Moreover, the fundamental group $G=\pi_1(M)$ 
has the structure of a semidirect product of free groups,  
$G=F_{n-1}\rtimes \Z$, and  $M$ is a $K(G,1)$. 
For details, see \cite{Z}, \cite{MS0}.

\begin{example} 
\label{ex=2134} 
Let $\A=\A(2134)$ be the arrangement defined in complex 
coordinates on $\RR^4=\C^2$ by the half-holomorphic function 
$Q(z,w)=zw(z-w)(z-2\bar{w})$; see Ziegler \cite[Example 2.2]{Z}. 
Using a computation from \cite[Example 5.10]{MS0}, we obtain 
the following presentation for the fundamental group of the complement
\[
G=\langle x_1,x_2,x_3,x_4 \mid 
(x_1,x_3^2x_4),\: (x_2,x_4),\: (x_3,x_4)\rangle. 
\]
It can be seen that  
$E_G=\widehat{\bL}(x_1,x_2,x_3,x_4)/
\llangle 2 [x_1,x_3] + [x_1,x_4], [x_2,x_4], [x_3,x_4] \rrangle$; 
thus, $G$ is $1$-formal.
The resonance variety $\R_1(G)\subset \C^4$ has two 
components, $\R^{\alpha}=\{ x \mid x_4=0\}$ and 
$\R^{\beta}=\{ x \mid x_4+2 x_3=0\}$, and the 
tangent cone formula holds for $G$.  
Though both components of $\R_1(G)$ are linear, the 
other three resonance obstructions are violated:
\begin{itemize}
\item The subspaces $\R^{\alpha}$ and  $\R^{\beta}$ 
are neither $0$-isotropic, nor $1$-isotropic.  
\item $\R^{\alpha}\cap \R^{\beta} = \{ x \mid x_3=x_4=0\}$, 
which is not equal to $\{0\}$. 
\item $\R_2(G)= \{x\mid x_1=x_3=x_4=0\} \cup 
\{x\mid x_2=x_3=x_4=0\}$, and neither of these components 
equals $\R^{\alpha}$ or $\R^{\beta}$.
\end{itemize}
Thus, $G$ is not the fundamental group of any smooth 
quasi-projective variety. 
\end{example}

Let $\A$ be an arrangement of transverse planes in $\RR^4$, 
with complement $M$. From the point of view of two classical 
invariants---the associated graded Lie algebra, and the 
Chen Lie algebra---the group $G=\pi_1(M)$ behaves like 
a $1$-formal group. Indeed,  the associated link $L$ 
has all linking numbers equal to $\pm 1$, in particular, 
the linking graph of $L$ is connected.  Thus, 
$\gr^*(G)\otimes \Q \cong \h_G$ and 
$\gr^*(G/G'')\otimes \Q \cong \h_G/\h''_G$, as graded 
Lie algebras, by \cite[Corollary 6.2]{MP} and 
\cite[Theorem 10.4(f)]{PS2}, respectively.  Nevertheless, 
our methods can detect non-formality, even in this 
delicate setting.

\begin{example} 
\label{ex=nonformal} 
Consider the arrangement $\A=\A(31425)$ defined 
in complex coordinates by the function 
$Q(z,w)=z(z-w)(z-2w)(2z+3w-5\overline{w}) (2z-w-5\overline{w})$; 
see \cite[Example 6.5]{MS1}. A computation shows that 
$TC_1(\V_2(G))$ has $9$ irreducible components, 
while $\R_2(G)$ has $10$ irreducible components; 
see \cite[Example 10.2]{MS2}, and \cite[Example 6.5]{MS1}, 
respectively.  By Theorem \ref{thm=tcf},  the group $G$ is 
not $1$-formal.  Thus, the complement $M$ cannot be a 
formal space, despite a claim to the contrary by 
Ziegler \cite[p.~10]{Z}.
\end{example}

\section{Configuration spaces} 
\label{sec=conf}

Denote by $S^{\times n}$ the $n$-fold cartesian product of a 
connected space $S$.  Consider the {\em configuration space}\/ 
of $n$ distinct labeled points in $S$,
\[
F(S, n)= S^{\times n}\setminus \bigcup_{i<j} \Delta_{ij} ,
\]
where $\Delta_{ij}$ is the diagonal $\{s\in S^{\times n} \mid  s_i=s_j\}$. 
The topology of configuration spaces has attracted considerable 
attention over the years.   For $S$ a smooth, complex projective 
variety, the cohomology algebra $H^*(F(S, n), \C)$ has been 
described by Totaro \cite{Tt}, solely in terms of $n$ and the 
cohomology algebra $H^*(S, \C)$. 

Let $C_g$ be a smooth compact complex curve of genus $g$ 
($g\ge 1$).  The fundamental group of the configuration space 
$M_{g,n}:= F(C_g, n)$ may be identified with $P_{g,n}$, the 
pure braid group on $n$ strings of the underlying Riemann surface.  
Starting from Totaro's description, it is straightforward to check 
that the low-degrees cup-product map of $P_{g,n}$ is equivalent, 
in the sense of Definition \ref{def=similar}, to the composite 
\begin{equation} 
\label{eq:cupg}
\mu_{g,n}\colon 
\xymatrix{ \bigwedge^2 H^1(C_g^{\times n},\C)
\ar[r]^(.55){\cup_{C_g^{\times n}}} & 
H^2(C_g^{\times n},\C)\ar@{>>}[r] &  
H^2(C_g^{\times n},\C)/ \spn \{[\Delta_{ij}]\}_{ i<j} } ,
\end{equation}
where $[\Delta_{ij}]\in H^2(C_g^{\times n}, \C)$ denotes the 
dual class of the diagonal $\Delta_{ij}$, and the second arrow 
is the canonical projection.  It follows that the connected  
smooth quasi-projective complex variety $M_{g,n}$ 
has the property that $W_1(H^1(M_{g,n},\C))=H^1(M_{g,n},\C)$, 
for all $g, n\ge 1$.

The Malcev Lie algebra of $P_{g,n}$ has been computed by 
Bezrukavnikov in \cite{B}, for all $g,n\ge 1$.  It turns out that 
the groups $P_{g,n}$ are $1$-formal, for $g>1$ and $n\ge 1$, 
or $g=1$ and $n\le 2$; see \cite[p.~130]{B}.   
On the other hand, Bezrukavnikov also states 
in \cite[Proposition 4.1(a)]{B} that $P_{1,n}$ is 
not $1$-formal for $n\ge 3$, without giving a full argument.  
With our methods, this can be easily proved.

\begin{example} 
\label{ex:g1not1f}
Let  $\{ a, b\}$ be the standard basis of $H^1(C_1,\C)=\C^2$. 
Note that the cohomology algebra 
$H^*( C_1^{\times n}, \C)$ is isomorphic to 
$\bigwedge^* (a_1, b_1,\dots ,a_n, b_n)$. 
Denote by $(x_1, y_1, \dots ,x_n ,y_n)$ the 
coordinates of $z\in H^1(P_{1,n}, \C)$. Using \eqref{eq:cupg}, 
it is readily seen that
\renewcommand{\arraystretch}{1.1}
\begin{equation*} 
\label{eq:resg1}
\R_1(P_{1,n})=\left\{ (x,y) \in \C^n\times \C^n \left|
\begin{array}{l}
\sum_{i=1}^n x_i=\sum_{i=1}^n y_i=0 ,\\
x_i y_j-x_j y_i=0,  \text{ for $1\le i<j< n$}
\end{array}
\right\}. \right.
\end{equation*}
\renewcommand{\arraystretch}{1.0}

Suppose $n\ge 3$.  Then $\R_1(P_{1,n})$ 
is a rational normal scroll in $\C^{2(n-1)}$, 
see \cite{Har}, \cite{E}. In particular, $\R_1(P_{1,n})$ is an 
irreducible, non-linear variety.  From Theorem \ref{thm=posob}, 
we conclude that $P_{1,n}$ is indeed non-$1$-formal. 
This indicates that Theorem 1.3 from \cite{Hai} 
cannot hold in the stated generality.

This family of examples also shows that both the 
$\R_1$--version of Arapura's result on $\V_1$ from 
Theorem \ref{thm:vadm}\eqref{va1} and the resonance 
obstruction test from our Theorem \ref{thm=posob}\eqref{a64} 
may fail, for an arbitrary smooth quasi-projective variety $M$. 
\end{example}

For $n\le 2$, things are even simpler. 

\begin{example} 
\label{ex:01iso}
It follows from \eqref{eq:cupg} that $\mu_{1,2}$ equals the 
canonical projection
\begin{equation*} 
\label{eq:cup12}
\mu_{1,2} \colon \bigwedge\nolimits^2 (a_1,b_1, a_2,b_2) \surj 
\bigwedge\nolimits^2 (a_1,b_1, a_2,b_2)/ \C \cdot 
(a_1-a_2)(b_1-b_2) .
\end{equation*}
It follows that $\R_1(P_{1,2})$ is a $2$-dimensional, 
$0$-isotropic linear subspace of $H^1 (P_{1,2},\C)$.

Consider now the smooth variety $M'_g:= M_{1,2}\times C_g$, 
with $g\ge 2$. By Proposition \ref{prop=products}, this 
variety has $1$-formal fundamental group. It also has 
the property that $W_1(H^1 (M'_g,\C))=H^1 (M'_g,\C)$. 
We infer from Lemma \ref{lem=resprod} that
\[
\R_1(\pi_1(M'_g))=\R_1(P_{1,2})\times\{0\} \cup
\{0\}\times  H^1 (C_g,\C) ,
\]
where the component $\R_1(P_{1,2})$ is $0$-isotropic 
and the component $H^1 (C_g,\C)$ is $1$-isotropic. 
We thus see that both cases listed in Proposition 
\ref{prop=posobstr}\eqref{ra1} may actually occur.
\end{example}

\begin{remark}
\label{rem:lib}
Recall from Example \ref{ex:notinj} that the tangent cone formula 
may fail for quasi-projective groups, at least in the case 
when $1$ is an isolated point of the characteristic variety.
The following statement can be extracted from \cite[p.~161]{Li}:
``If $M$ is a quasi-projective variety and $1$ is not an isolated 
point of $\V_1(\pi_1(M))$, then $TC_1(\V_1(\pi_1(M)))= 
\R_1(\pi_1(M))$."   Taking $M$ to be one of the configuration 
spaces $M_{1,n}$, with  $n\ge 3$, shows that this statement 
does not hold, even when $M$ is smooth.

Indeed, since $P_{1,2}$ is $1$-formal, we obtain from 
Theorem \ref{thm=tcfintro} that $\V_1(P_{1,2})$ is 
$2$-dimensional at $1$. As is well-known,
the natural surjection, $P_{1,n}\surj P_{1,2}$, embeds
$\V_1(P_{1,2})$ into $\V_1(P_{1,n})$, for $n\ge 2$. Thus, 
$\V_1(P_{1,n})$ is positive-dimensional at $1$, for $n\ge 2$.
On the other hand, it follows from Example \ref{ex:g1not1f} that
$TC_1(\V_1(P_{1,n}))$ is strictly contained in $\R_1(P_{1,n})$, 
for $n\ge 3$.
\end{remark}

\section{Artin groups} 
\label{sec=artinalg}

In this section, we analyze the class of finite-type Artin 
groups. Using the resonance obstructions from 
Theorem \ref{thm=posobstr}, we give a complete 
answer to Serre's question for right-angled Artin groups, 
and we give a Malcev Lie algebra version of the answer for 
arbitrary Artin groups. 

\subsection{Labeled graphs and Artin groups}
\label{subsec=artin}
Let $\Gamma=(\sV,\sE,\ell)$ be a labeled finite simplicial graph, 
with vertex set $\sV$, edge set $\sE\subset \binom{\sV}{2}$, 
and labeling function $\ell \colon \sE\to \N _{\geq 2}$.  
Finite simplicial graphs are identified in the sequel with labeled 
finite simplicial graphs with $\ell (e)=2$, for each $e \in\sE$.

\begin{definition} 
\label{def=artgps} 
The {\em Artin group}\/ associated to the labeled graph $\Gamma$ is 
the group $G_{\Gamma}$ generated by the vertices $v \in \sV$ and 
with a defining relation 
\[
\underbrace{vwv\cdots}_{\ell(e)} =\underbrace{wvw\cdots}_{\ell(e)}
\]
for each edge $e=\{v,w\}$ in $\sE$.  If $\Gamma$ is unlabeled, 
then $G_{\Gamma}$ is called a {\it right-angled Artin group}, and 
is defined by commutation relations $vw=wv$,
one for each edge $\{v,w\} \in\sE$.
\end{definition}

\begin{example}
\label{ex=artjoin}
Let $\Gamma=(\sV,\sE,\ell)$ and $\Gamma'=(\sV',\sE',\ell')$ be two 
labeled graphs.  Denote by $\Gamma \bigsqcup \Gamma '$ their 
disjoint union, and by $\Gamma *\Gamma '$ their {\em join}, 
with vertex set $\sV \bigsqcup \sV '$, edge set 
$ \sE \bigsqcup \sE ' \bigsqcup \{  \{ v,v'\} \mid  v \in \sV, v' \in \sV'\}$, 
and label $2$ on each edge $ \{ v,v'\}$.  Then
\[
G_{\Gamma \bigsqcup \Gamma '}  = G_{\Gamma} *G_{\Gamma '}
\quad \text{and} \quad 
G_{\Gamma *\Gamma '}  = G_{\Gamma} \times G_{\Gamma '}\, .
\]
In particular, if $\Gamma$ is a discrete graph, i.e., $\sE=\emptyset$, 
then $G_{\Gamma}=\FF_n$, whereas if $\Gamma$ is an (unlabeled) 
complete graph, i.e., $\sE=\binom{\sV}{2}$, then 
$G_{\Gamma}=\Z^n$, where $n=\abs{\sV}$.  More generally, 
if $\Gamma$ is a complete multipartite graph  (i.e., a finite 
join of discrete graphs), then $G_{\Gamma}$ is a finite 
direct product of finitely generated free groups. 
\end{example}

Given a graph $\Gamma=(\sV,\sE)$ and a subset of vertices $\sW \subset \sV$, 
we denote by $\Gamma(\sW)$ the full subgraph of $\Gamma$, 
with vertex set  $\sW $ and edge set $\sE \cap \binom{\sW}{2} $.

Let $(S^1)^{\sV}$ be the compact $n$-torus, where $n=|\sV|$, endowed with 
the standard cell structure.  Denote by $K_{\Gamma}$ the subcomplex of 
$(S^1)^{\sV}$ having a $k$-cell for each subset $\sW \subset \sV$ of size $k$ 
for which $\Gamma(\sW)$ is a complete graph. As shown by Charney--Davis
\cite{CD} and Meier--VanWyk \cite{MV},  $K_{\Gamma}=K(G_{\Gamma},1)$. 
In particular, the cup-product map 
$\cup_{G_{\Gamma}}\colon H^1(G_{\Gamma},\C)\wedge H^1(G_{\Gamma},\C) 
\to H^2(G_{\Gamma},\C)$ may be identified with the linear map 
$\cup_{\Gamma}\colon \C^\sV \wedge  \C^\sV \to \C^\sE$ 
defined by
\begin{equation} 
\label{eq=cupart}
v\cup_{\Gamma}w = 
\begin{cases} 
\pm  \{v,w\} , &\text{if  $\{v,w\} \in \sE$},\\
0,  &\text{otherwise} ,
\end{cases}  
\end{equation}
with signs determined by fixing an orientation of the edges of $\Gamma$.

\subsection{Jumping loci for right-angled Artin groups}
\label{subsec=jumpartin}

The resonance varieties of right-angled Artin groups  
were  described explicitly in Theorem 5.5 from \cite{PS1}. 
If $\Gamma=(\sV,\sE)$ is a graph, then 
\begin{equation} 
\label{eq=resart}
\R_1(G_{\Gamma})= \bigcup_{\sW} \C^ \sW ,
\end{equation}
where the union is taken over all subsets $\sW \subset \sV$ 
such that $\Gamma(\sW)$ is disconnected, and maximal with 
respect to this property.   Moreover, the decomposition \eqref{eq=resart} 
coincides with the decomposition into irreducible components of 
$\R_1(G_{\Gamma})$.

Before proceeding to the Serre problem, we describe the 
characteristic variety of $G_{\Gamma}$. 
For $\sW \subset \sV$, define the subtorus 
$\T _{\sW} \subset \T _{G_{\Gamma }   }=(\C^*)^{\sV}$ by
\[
\T _{\sW}=\{(t_v)_{v \in \sV} 
\in (\C^*)^{\sV} ~ | ~ t_v=1 ~ \text{for} ~ v \notin \sW \} .
\]
The  map 
$\exp\colon T_1\T_{G_\Gamma} \to \T_{G_{\Gamma}}$ 
is the componentwise  exponential map 
$\exp^{\sV}\colon \C^\sV\to(\C^*)^{\sV}$; 
its restriction to the subspace spanned by $\sW$ 
is $\exp^{\sW}\colon \C^\sW\to(\C^*)^{\sW}=\T_{\sW}$.

\begin{prop} 
\label{prop=vart}
Let $G_{\Gamma }$ be the right-angled Artin group associated 
to the graph $\Gamma=(\sV,\sE)$. Then
\[
\V_1(G_{\Gamma })= \bigcup _{\sW}\T _{\sW} ,
\]
where the union is over all subsets $\sW\subset \sV$ 
such that $\Gamma(\sW)$ is maximally disconnected.   
Moreover, this decomposition coincides with the 
decomposition into irreducible components of 
$\V_1(G_{\Gamma})$.

\end{prop}

\begin{proof}
The realization of $K(G_{\Gamma },1)$ as a subcomplex $K_{\Gamma }$ 
of the torus $(S^1)^{\sV}$ yields a well-known resolution of the trivial 
$\Z G_{\Gamma}$-module $\Z$ by finitely generated, free 
$\Z G_{\Gamma}$-modules, as the augmented, 
$G_{\Gamma}$-equivariant chain complex of the 
universal cover of $K_{\Gamma }$,
\[
{\widetilde C}_{\bullet} (\widetilde{K_{\Gamma}}) \colon 
\xymatrixcolsep{18pt}
\xymatrix{\cdots \ar[r] &  \Z G_{\Gamma } \otimes C_k  
\ar^(.45){d_k }[r]& \Z G_{\Gamma } \otimes C_{k-1}  
\ar[r] & \cdots \ar^(.45){d_1}[r] & \Z G_{\Gamma } 
\ar^{\epsilon}[r] & \Z \ar[r] & 0  }. 
\]
Here $C_k$ denotes the free abelian group generated by the 
$k$-cells of $K_{\Gamma }$, and the boundary maps are given by 
\begin{equation} 
\label{eq=kchains}
d_k(e_{v_1} \times \cdots \times e_{v_k})=\sum _{i=1}^{k}(-1)^{i-1}(v_i-1) \otimes 
e_{v_1} \times \cdots \times {\widehat  e_{v_i} } \times \cdots \times e_{v_k} ,
\end{equation}
where, for each $v \in \sV$, the symbol $e_v$ denotes   
the $1$-cell  corresponding to the $v$-th factor of 
$(S^1)^{\sV}$.

Now let $\rho= (t_v)_{v \in \sV} \in (\C^*)^{\sV}$ be an arbitrary character. 
Denoting by $\{v^*\} _{v \in \sV}$ the basis of $H^1(G_{\Gamma },\C)$ 
dual to the canonical basis of $H_1(G_{\Gamma },\C)$, define 
an element $z \in \C^\sV= H^1(G_{\Gamma },\C)$ by
$z=\sum_{v \in \sV}(t_v-1)v^*$.  Using \eqref{eq=kchains}, 
it is not difficult to check the following equality of cochain 
complexes
\begin{equation} 
\label{eq=reschar}
\Hom_{\Z  G_{\Gamma }}( {\widetilde C}_{\bullet} 
(\widetilde{K_{\Gamma}}),{}_{\rho}\C) = 
(H^{\bullet} (G_{\Gamma },\C), \mu_ z ) .
\end{equation}
It follows then, directly from the definitions \eqref{eq=charm} 
and \eqref{eq=resmu}, and using \eqref{eq=reschar}, that 
$\rho \in \V_1(G_{\Gamma })$ if and only if 
$z \in \R_1(G_{\Gamma })$.  Hence, the claimed 
decomposition of $\V_1(G_{\Gamma })$ is a direct 
consequence of the decomposition \eqref{eq=resart}.
\end{proof}

\subsection{Serre's problem for right-angled Artin groups}
\label{subsec=rightartin}

As shown by Kapovich and Millson in \cite[Theorem 16.10]{KM}, 
all Artin groups are $1$-formal.  This opens the way for approaching 
Serre's problem for Artin groups via resonance varieties, which, 
as noted above, were described explicitly in \cite{PS1}.  Using 
these tools, we find out precisely which right-angled Artin groups 
can be realized as fundamental groups of quasi-compact 
K\"{a}hler manifolds.

\begin{theorem} 
\label{thm=raserre} 
Let  $\Gamma=(\sV,\sE)$ be a finite simplicial graph, with associated 
right-angled Artin group $G_{\Gamma}$. The following are equivalent.

\begin{romenum}
\item  \label{s1} 
The group $G_{\Gamma}$ is quasi-K\"{a}hler.
\item  \label{s2} 
The isotropicity property from Definition \ref{def=posob}
is satisfied by $\cup_{G_{\Gamma}}$.
\item  \label{s3} 
The graph  $\Gamma$ is complete multipartite graph.
\item  \label{s4}  
The group $G_{\Gamma}$ is a product of finitely generated 
free groups.
\end{romenum}
\end{theorem}

\begin{proof}
For the implication \eqref{s1} $\Rightarrow$ \eqref{s2}, use 
the $1$-formality of $G_{\Gamma}$ and 
Theorem \ref{thm=posobstr}.

\smallskip
The implication  \eqref{s2} $\Rightarrow$ \eqref{s3} 
is proved by induction on $n=\abs{\sV}$. If  $\Gamma$ is complete, 
then  $\Gamma$ is the join of $n$ graphs with one vertex. Otherwise, 
there is a subset $\sW \subset \sV$ such that $\Gamma (\sW)$ is 
disconnected, and maximal with respect to this property.  Write 
$\sW = \sW' \bigsqcup   \sW''$, with both $ \sW'$ and $ \sW''$ 
non-empty and with no edge connecting $ \sW'$ to $ \sW''$. 
Then $\Gamma (\sW)=\Gamma (\sW ')  \bigsqcup \Gamma (\sW '')$, 
and so $G_{\Gamma (\sW)}=G_{\Gamma (\sW ')} * 
G_{\Gamma (\sW '')} $.  Hence, $w' \cup_{\Gamma (\sW)} w''=0$, 
for any $w' \in \sW '$ and $w'' \in \sW''$.  We infer from  
\cite[Lemma 5.2]{PS1} that $\R_1(G_{\Gamma (\sW)})=\C^\sW$. 

On the other hand, we know from \eqref{eq=resart} that $\C^\sW$
is a positive-dimensional irreducible component of 
$\R_1(G_{\Gamma })$. Our hypothesis implies that 
$\C^\sW$ is either $0$-isotropic  or $1$-isotropic with 
respect to $\cup_{\Gamma (\sW)}$.  By Lemma \ref{lem=reswedge}, 
$ \cup _{\Gamma (\sW ')}=\cup _{\Gamma (\sW '')}  =0$.  
The cup-product formula \eqref{eq=cupart} implies that 
$\Gamma (\sW)$ is a discrete graph.

If $\sW=\sV$, we are done. Otherwise, $\sV=\sW \bigsqcup 
\sW_1$, with $\abs{\sW_1} <n$. Since  $\Gamma (\sW)$ is maximally 
disconnected, this forces $\{ w, w_1\} \in \sE$, for all $w \in \sW$ 
and $w_1 \in \sW_1$. In other words, $\Gamma $ is the join 
$\Gamma (\sW) * \Gamma (\sW_1)$; thus, 
$G_{\Gamma }=G_{\Gamma (\sW)} \times G_{\Gamma (\sW_1)}$. 
By Lemma \ref{lem=resprod}, $\cup_{\Gamma (\sW_1)}$ inherits 
from $\cup_{\Gamma}$ the isotropicity property.  
This completes the induction.

\smallskip
The implication  \eqref{s3} $\Rightarrow$ \eqref{s4} follows 
from the discussion in Example \ref{ex=artjoin}.

\smallskip
Finally, the implication  \eqref{s4} $\Rightarrow$ \eqref{s1} 
follows by taking products and realizing free groups by the 
complex line with a number of points deleted.
\end{proof}

As is well-known, two right-angled Artin groups are isomorphic 
if and only if the corresponding graphs are isomorphic. 
Evidently, there are infinitely many graphs which are not 
joins of discrete graphs.  Thus, implication \eqref{s1}  
$\Rightarrow$ \eqref{s3} from Theorem \ref{thm=raserre} 
allows us to recover, in sharper form, a result of Kapovich 
and Millson (Theorem 14.7 from \cite{KM}).

\begin{corollary}
\label{cor=rkm}
Among right-angled Artin groups $G_{\Gamma}$, there are infinitely 
many mutually non-isomorphic groups which are not isomorphic 
to fundamental groups of smooth, quasi-projective complex varieties. 
\end{corollary}

\subsection{A Malcev Lie algebra version of Serre's question}
\label{subsec=mlserre}

Next, we describe a construction that associates to a labeled 
graph $\Gamma=(\sV,\sE,\ell)$ an ordinary graph, 
${\tilde \Gamma}=({\tilde  \sV}, {\tilde \sE})$, which 
we call the {\em odd contraction}\/ of $\Gamma$.
First define $\Gamma_{\odd}$ to be the unlabeled 
graph with vertex set $\sV$ and edge set 
$\{e \in \sE \mid  \ell (e)~  \text{is odd}\}$.
Then define ${\tilde  \sV}$ to be the set of connected 
components of $\Gamma_{\odd}$, with two distinct 
components determining an edge $\{c,c'\} \in  {\tilde \sE}$ 
if and only if there exist vertices $v\in c$ 
and $v'\in c'$ which are connected by an edge in $\sE$.

\begin{example} 
\label{ex=braids} 
Let $\Gamma$ be the complete graph on vertices $\{1,2,\dots ,n-1\}$, 
with labels $\ell(\{i,j\})=2$ if $|i-j|>1$ and  
$\ell(\{i,j\})=3$ if $|i-j|=1$. The corresponding Artin group 
is the classical  braid group on $n$ strings, $B_n$. 
Since in this case  $\Gamma_{\odd}$ is connected, the odd 
contraction ${\tilde \Gamma}$ is the discrete graph with 
a single vertex.
\end{example}

\begin{lemma}
\label{lem=malcont}
Let $\Gamma=(\sV,\sE,\ell)$ be  a labeled graph, with odd 
contraction ${\tilde \Gamma}=({\tilde  \sV} , {\tilde \sE})$. 
Then the Malcev Lie algebra of $G_{\Gamma}$ is filtered Lie 
isomorphic to the Malcev Lie algebra of  $G_{\tilde \Gamma}$.
\end{lemma}

\begin{proof}
The Malcev Lie algebra of $G_{\Gamma }$ was computed in  
\cite[Theorem 16.10]{KM}. It is the quotient of the free 
Malcev Lie algebra on $\sV$, ${ \widehat \bL}(\sV)$, by 
the closed Lie ideal generated by the differences $u-v$, 
for odd-labeled edges $\{u,v\}\in \sE$, and by the 
brackets $[u,v]$, for even-labeled edges $\{u,v\}\in \sE$. 
Plainly, this quotient is filtered Lie isomorphic to the 
quotient of  ${ \widehat \bL}({\tilde \sV})$ by the closed 
Lie ideal generated by the brackets $[c,c']$, for 
$\{c,c'\} \in {\tilde \sE}$, which is just the Malcev Lie 
algebra of $G_{ {\tilde \Gamma}    }$. 
\end{proof}

The Coxeter group associated to a labeled graph 
$\Gamma=(\sV,\sE,\ell)$ is the quotient of the Artin group 
$G_{\Gamma}$ by the normal subgroup generated by 
$\{v^2 \mid v\in \sV\}$.  If the Coxeter group $W_{\Gamma}$ 
is finite, then $G_{\Gamma}$ is quasi-projective. 
The proof of this assertion, due to Brieskorn \cite{Br}, 
involves the  following steps:  $W_\Gamma$ acts as a 
group of reflections in some $\C^{n}$; the action is free 
on the complement $M_{\Gamma}$ of the arrangement 
of reflecting hyperplanes of $W_{\Gamma}$, and 
$G_\Gamma=\pi_1(M_{\Gamma}/W_\Gamma)$;   
finally, the quotient manifold $M_{\Gamma}/W_\Gamma$ 
is a complex smooth quasi-projective variety.

It would be interesting to know exactly which (non-right-angled) 
Artin groups can be realized by smooth, quasi-projective complex 
varieties.  We give an answer to this question, albeit only at the 
level of Malcev Lie algebras of the respective groups.

\begin{corollary} 
\label{cor=aserre}
Let $\Gamma$ be  a labeled graph, with associated 
Artin group $G_{\Gamma }$ and odd contraction the 
unlabeled graph ${\tilde \Gamma}$. The following are equivalent.

\begin{romenum}
\item \label{arba1}
The Malcev Lie algebra of $G_{\Gamma}$ is filtered Lie 
isomorphic to the Malcev Lie algebra of a quasi-K\"{a}hler group.

\item \label{arba2}
The isotropicity property from Definition \ref{def=posob}
is satisfied by $\cup_{G_{\Gamma}}$.

\item  \label{arba3} 
The graph ${\tilde \Gamma}$ is a complete multipartite graph.

\item  \label{arba4}  
The Malcev Lie algebra of $G_{\Gamma}$ is filtered Lie 
isomorphic to the Malcev Lie algebra of a product of finitely 
generated free groups.
\end{romenum}
\end{corollary}

\begin{proof}
By Lemma \ref{lem=malcont}, the Malcev Lie algebras of 
$G_{\Gamma}$ and $G_{ {\tilde \Gamma} }$ are filtered isomorphic. 
Hence, the graded Lie algebras $\gr^*(G_{\Gamma }) \otimes \C$ 
and  $\gr^*(G_{ {\tilde \Gamma} }) \otimes \C$ are isomorphic.

From \cite{S}, we know that the kernel of 
the Lie bracket, $\bigwedge ^2 \gr^1(G) \otimes \C \to 
\gr^2(G) \otimes \C$, is equal to  $\im (\partial _G)$, 
for any finitely presented group $G$. It follows that the 
cup-product maps $\cup _{G_{\Gamma }}$ and 
$\cup _{G_{ {\tilde \Gamma}}}$ are equivalent, in the 
sense of Definition \ref{def=similar}.  Consequently, 
$\cup_{G_{\Gamma}}$ satisfies the second resonance 
obstruction if and only if $\cup_{G_{{\tilde\Gamma}}}$ 
does so.

With these remarks, the Corollary follows at once 
from Theorems \ref{thm=raserre} and \ref{thm=posob}.
\end{proof}

\subsection{K\"{a}hler right-angled Artin groups}
\label{subsec:kraag}

With our methods, we may easily decide which 
right-angled Artin groups are K\"{a}hler groups.  

\begin{corollary} 
\label{cor:raag kahler} 
For a right-angled Artin group $G_{\Gamma}$, 
the following are equivalent.
\begin{romenum}
\item  \label{agk1} 
The group $G_{\Gamma}$ is K\"{a}hler.
\item  \label{agk2} 
The graph  $\Gamma$ is a complete graph on an even 
number of vertices.
\item  \label{agk3}  
The group $G_{\Gamma}$ is a free abelian group of 
even rank. 
\end{romenum}
\end{corollary}

\begin{proof}
Implications \eqref{agk2} $\Rightarrow$ \eqref{agk3} 
$\Rightarrow$ \eqref{agk1}  are clear. So suppose 
$G_{\Gamma}$ is a K\"{a}hler group. By Theorem 
\ref{thm=raserre}, $\Gamma$ is a complete multi-partite graph 
$\overline{K}_{n_1}*\cdots * \overline{K}_{n_r}$, and 
$G_\Gamma=\bF_{n_1}\times \cdots \times \bF_{n_r}$.  
By Lemma \ref{lem=resprod}, and abusing notation 
slightly, $\R_1(G_{\Gamma })=\bigcup_i \R_1(\bF_{n_i})$.
Now, if there were an index $i$ for which $n_i>1$, then 
$\R_1(\bF_{n_i})=\C^{n_i}$ would be a positive-dimensional, 
$0$-isotropic, irreducible component of $\R_1(G_{\Gamma })$, 
contradicting Corollary \ref{cor:firstex}\eqref{appl1}. 
Thus, we must have $n_1=\cdots =n_r=1$, and $\Gamma=K_r$.  
Moreover, since $G_\Gamma=\Z^r$ is a K\"{a}hler group, 
$r$ must be even. 
\end{proof}

\begin{ack}
The second and third authors are grateful to the
Mathematical Sciences Research Institute in Berkeley, California,
for its support and hospitality during the program ``Hyperplane
Arrangements and Applications'' (Fall, 2004).
A preliminary version of this paper was presented 
by the second author in a talk delivered at the MSRI
workshop on ``Topology of Arrangements and Applications''
(October 4--8, 2004). Part of this work was done during the 
second author's visit at Universit\'{e} de Nice--Sophia Antipolis 
in Nice, France, which provided an excellent environment 
for research (April--June, 2005).
\end{ack}

\quad\newline\vspace{-1.1in}
\newcommand{\arxiv}[1]
{\texttt{\href{http://arxiv.org/abs/#1}{arXiv:#1}}}

\renewcommand{\MR}[1]
{\href{http://www.ams.org/mathscinet-getitem?mr=#1}{MR#1}}

\bibliographystyle{amsplain}

\end{document}